\documentclass[11pt]{amsart}

\usepackage{amsfonts}
\usepackage{amssymb}
\usepackage{graphicx}
\usepackage{color}

\def\nn{\mathbb{N}}
\def\rr{\mathbb{R}}

\def\C{\mathcal{C}}
\def\D{\mathcal{D}}
\def\G{\mathcal{G}}

\def\F{\mathcal{F}}
\def\I{\mathcal{I}}
\def\J{\mathcal{J}}
\def\L{\mathcal{L}}

\def\T{\mathcal{T}}
\def\bI{\bar{\mathcal I}}
\def\bJ{\bar{\mathcal J}}
\def\GG{{\frak G}}

\def\epsilon{\varepsilon}
\def\leq{\leqslant}
\def\geq{\geqslant}

\def\x{{\bf x}}
\def\y{{\bf y}}
\def\k{{\bf k}}

\def\pp{\mathbb{P}}

\def\ms{{\medskip\noindent}}
\def\bs{{\bigskip\noindent}}

\def\per{{\rm Per}}
\def\bx{\bar{\x}}
\def\tx{\tilde{\x}}
\def\by{\bar{\y}}
\def\ty{\tilde{\y}}
\def\bT{\bar{T}}
\def\tT{\tilde{T}}
\def\bsigma{\bar{\sigma}}
\def\bG{\bar{\G}}
\def\bV{\bar{V}}
\def\bA{\bar{A}}
\def\bu{\bar{u}}

\def\Id{{\rm Id}}
\def\osc{{\rm osc}}
\def\Idosc{{\rm Id_{osc}}}
\def\dist{{\rm d}_{\G}}
\def\lcm{{\rm lcm}}
\def\vol{{\rm vol}}
\def\ext{{\rm ext}}

\newtheorem{proposicion}{Proposition}

\newtheorem{teorema}{Theorem}
\newtheorem{lema}{Lemma}
\newtheorem{nota}{Remark}
\newtheorem{corolario}{Corollary}
\def\endproof{\hfill$\Box$}

\begin{document}

\title[Regulatory Dynamics on Random Networks]{Regulatory Dynamics on Random Networks: Asymptotic Periodicity and Modularity}
\author{A.~Cros, A.~Morante, and E.~Ugalde}

\maketitle

\begin{center}

\begin{minipage}{14cm}
\begin{small}
CUCEI, Universidad de Guadalajara, 44430 Guadalajara, M\'exico,\\
F.~Ciencias, Universidad Aut\'onoma de San Luis Potos\'\i , 
78000 San Luis Potos\'\i, M\'exico, \& \\
I.~F\'\i sica, Universidad Aut\'onoma de San Luis Potos\'\i , 
78000 San Luis Potos\'\i, M\'exico.
\end{small}
\end{minipage}

\vspace{8truemm} 

\begin{abstract} We study the dynamics of discrete--time regulatory 
networks on random digraphs. For this we define ensembles of 
deterministic orbits of random regulatory networks, and introduce 
some statistical indicators related to the long--term dynamics of the system.
We prove that, in a random regulatory network, initial conditions converge
almost surely to a periodic attractor. 
We study the subnetworks, which we call modules, where the 
periodic asymptotic oscillations are concentrated. 
We proof that those modules are dynamically equivalent to independent 
regulatory networks. 

\vskip 2cm 
\ms \noindent{\it Keywords}: regulatory networks, random graphs, 
coupled map networks.

\end{abstract}
\end{center}

\medskip

\bs\section{Introduction}

\noindent Numerous natural and artificial systems can be though as a 
collection of basic units interacting according to simple rules. 
Examples of this interacting systems are the genetic regulatory 
networks, composed of interactions between DNA, RNA, proteins, and 
small molecules. In social or ecological networks, a similar regulatory 
dynamics may also be considered. The traditional way to model these 
systems is by using coupled differential equations, and more 
particularly systems of piecewise affine differential 
equations~(see~\cite{deJong02, Edwards00, Mestl&al95}). Finite state models, 
better known as logical networks, are also used 
(see~\cite{Kauffman69, Glass&Kauffman74, Thieffry&Thomas95}). Within 
these modeling strategies, the interacting units have a regular behavior 
when taken separately, but are capable to generate global complex 
dynamics when arranged in a complex interaction architecture. 
In all the models considered so far, each interacting unit regulates 
some other units in the collection by enhancing or repressing their 
activity. It is possible then to define an underlying network with 
interacting units as vertices, and their interactions as 
arrows connecting those vertices. The theoretical problem we face here 
is to understand the relation between the structure of the underlying 
network, and the possible dynamical behaviors of the system. We will do
this in the context of a particular class of models first 
introduced in~\cite{Volchenkov&Lima05}, and further studied in~\cite{Coutinho&al06} 
and~\cite{Lima&Ugalde06}. In these models, the level of activity of each 
unit is codified by a positive real number. The system evolves synchronously 
at discrete time steps, each unit following an affine contraction dictated by 
the activity level and interaction mode of its neighboring units. 
The contraction coefficient of those transformations determines the 
degradation rate at which, in absence of interactions, the activity of a 
given unit vanishes. In the framework of this modeling we have proved  
general results concerning the constrains imposed by the structure of 
the underlying network, over the possible asymptotic behaviors of a fixed 
system~\cite{Lima&Ugalde06}. In the present paper, following~\cite{Volchenkov&Lima05}, 
we will focus on the asymptotic dynamics of regulatory systems whose interactions 
are chosen at random at the beginning of the evolution. Within this approach, 
individual orbits are elements of a sample space, and the statistical 
indicators we will study become orbit dependent random variables. The probability 
measures we used are built from a fixed probability distribution over the set of
possible underlying networks. Then, given a fixed underlying network, we 
associate a sign to each one of its arrows, depending on whether the interaction 
they represent are activations or inhibitions. Positive and negative signs are 
randomly chosen, keeping a fixed proportion of negative arrows inside a given 
statistical ensemble of systems. In this way it is possible to study 
certain characteristics of the asymptotic dynamics, as function of the degradation
rate and the proportion of inhibitory interactions. 

\ms Our first result states that in a random regulatory network, 
initial conditions converge almost surely to a periodic attractor. This result
points to the conclusion that in regulatory dynamics, the origin of the complexity
is the coexistence of multiple dynamically simple attractors. We prove that this is
the case in a full measure set in the parameter space. 

\ms According to our preliminary numerical explorations, the long--term 
oscillations of the system concentrate on subnetwork whose structure 
depends on the parameters of the statistical ensemble of regulatory networks. 
Our second result states that the dynamics one can observe when 
restricted to the oscillatory subnetwork, is equivalent 
to the dynamics supported by the subnetwork considered as an isolated system.  
This result allows us to introduce the concept of modularity. If we call module 
any observable oscillatory subnetwork, then, according to our result, the dynamics of a 
small network is preserved when it appears as a module in a larger network. 
We interpret this as the emergence of modularity.
This result allows us to predict admissible asymptotic behaviors in regulatory 
networks admitting disconnected oscillatory subnetworks. 
It is worth mentioning that this kind of modularity was already studied in 
the context of Boolean networks~\cite{Bastolla&Parisi98} and more recently,
in continuous--time regulatory networks~\cite{Gomez&al06}. Our approach allows 
a formal approach to the problems addressed in those works.

%%%%%%%%%%%  ORGANIZACION DEL RESTO DEL ARTICULO
\ms The paper is organized as follows. In the next section we will introduce the
objects under study, then in Section~\ref{section-results} we will state the 
results, and present some of the proofs. After reviewing two examples, which we do 
in Section~\ref{section-examples}, we will give the proofs of the two more technical
results in Section~\ref{section-proofs}. The paper ends with a section of final 
comments and conclusions.

%\ms The paper is organized as follows. In the next section we will introduce the
%objects under study, then in Section~\ref{section-results} we will state the 
%results, and present some of the proofs. After reviewing two examples, which we 
%do in Section~\ref{section-examples}, we end the paper with a section of final comments and conclusions.
%In Appendix~\ref{section-appendix} we will give the proofs of the two more technical results.

\ms This work was supported by CONACyT through the grant SEP--2003--C02--42765,
and by the cooperation ECOS--CONACyT--ANUIES M04--M01. E.~U.~thanks Bastien Fernandez
for his suggestions and comments. 

\bs\section{Preliminaries}~\label{section-preliminaries}

\ms \subsection{Regulatory Networks as Dynamical 
Systems}~\label{subsection-regulatory-networks}\

\ms The interaction architecture of the regulatory network is encoded in 
a directed graph $\G=(V, A)$, where vertices $V$ represent interacting units, 
and the arrows $A\subset V\times V$ denote interaction between them. 
To each interaction $(u,v)\in A$ we associate a threshold $T_{(u,v)} \in [0,1]$, 
and a sign $\sigma_{(u,v)}\in\{-1,1\}$ that is chosen according to whether this 
interaction is an inhibition or an activation.  
We quantify the activity of each unit $v\in V$ with real number $\x_v\in [0,1]$. 
Thus, the activation state of the network at a given time $t$ is determined by 
the vector $\x^{t}\in [0,1]^V$. The influence of a unit $u$ over a target 
unit $v$ turns on or off, depending on its sign, when the value of $\x_u$ 
trespasses the threshold $T_{(u,v)}$. The evolution of the network is generated 
by the iteration of the map $F_{\G,\sigma,T,a}:[0,1]^V\to[0,1]^V$ such that 
\begin{equation}\label{transformation}
\x^{t+1}=F_{\G,\sigma,T,a}(\x^t):=a\x^t+(1-a)D_{\G,\sigma,T}(\x^t), 
\end{equation}  
where the contraction rate $a\in [0,1)$ determines the speed of degradation of the 
activity of the units in absence of interaction, and the interaction term 
$D_{\G,\sigma,T}:[0,1]^V\to [0,1]^V$ is the piecewise constant function defined by
\begin{equation}\label{constant-part}
D_{\G,\sigma,T}(\x)_v:=\frac{1}{\Id(v)}\sum_{u\in V: (u,v)\in A} 
H(\sigma_{(u,v)}\left(\x_u-T_{(u,v)})\right).
\end{equation}
Here $H:\rr\to\rr$ is the Heaviside function, and $\Id(v):=\#\{u\in V: (u,v)\in A\}$ 
stands for the input degree of the vertex $v$. We will be referring to the 
discontinuity set of the transformation $D_{\G,\sigma,T}$, which is 
\begin{equation}\label{discontinuity-set}
\Delta_T:=\left\{\x\in[0,1]^V:\ \x_u=T_{u,v}\ \text{for some} (u,v)\in A\right\}.
\end{equation}
For each $a\in [0,1)$, the transformation $F_{\G,\sigma,T,a}$ is a piecewise 
affine contraction with discontinuity set $\Delta_T$. 

\ms From now on, by a discrete--time regulatory network we will mean a 
discrete--time dynamical system $([0,1]^V,F_{\G,\sigma,T,a})$, with phase 
space $[0,1]^V$, and evolution generated by the piecewise affine contraction 
$F_{\G,\sigma,T,a}:[0,1]^V\to [0,1]^V$ defined in Equation~(\ref{transformation}). 
The discrete--time regulatory networks studied here have interactions of equal
strength, and they act additively on each target unit. More general discrete--time 
regulatory networks have been considered in~\cite{Coutinho&al06, Lima&Ugalde06}.

\bs\subsection{Statistical Ensembles}~\label{subsection-statistical-ensemble}\

\ms We build our statistical ensembles as follows. We fix the 
value of the contraction rate $a\in[0,1)$ and the set $V$ representing the 
interacting units. Then we consider the set of all possible piecewise transformations,
\begin{equation}\label{equation-transformations-set}
\F_{a,V}:=\left\{F_{\G,\sigma,T,a}:\ \G:=(V,A), \ \sigma\in \{-1,1\}^A, \ T\in[0,1]^A\right\}.
\end{equation}
The individual elements of our statistical ensembles are couples 
$(F_{\G,\sigma,T,a},\x)$, with $F_{\G,\sigma,T,a}\in \F_{a,V}$ and $\x\in [0,1]^V$. 
A couple $(F_{\G,\sigma,T,a},\x)$ determines a deterministic orbit 
$\{\x^t:=F_{\G,\sigma,T,a}^t(\x)\}_{t=0}^{\infty}$. 
 
\ms We supply the sample space $\F_{a,V}\times [0,1]^V$ with a probability measure 
as follows. First we fix a probability distribution $\pp_G$ over the set 
$\GG_V:=\{(V, A): A\subset V\times V\}$ of all directed graphs with vertex set $V$.  
Then, for $\eta\in [0,1]$ and $\G=(V,A)$, we choose sign $-1$ with probability $\eta$ 
and $+1$ with probability $1-\eta$, independently for each arrow in $\G$. The 
thresholds are independent and uniformly distributed random variables in $[0,1]^A$, 
as well as the initial conditions in $[0,1]^V$. In this way we obtain the probability 
measure $\pp_{a,\eta}$ on $\F_{a,V}\times [0,1]^V$ such that 
\begin{equation}\label{equation-measure}
\pp_{a,\eta}\left\{(F_{\G,\sigma,T,a},\x):\ T\in \I, \x\in \J\right\}
=\pp_G(\G)\times \pp_{A,\eta}(\sigma)\times \vol(\I)\times \vol(\J),
\end{equation}
for all rectangles $\I\subset [0,1]^A$ and $\J\subset [0,1]^V$. Here 
$\pp_{A,\eta}:\{-1,1\}^A\to [0,1]$ is such that 
\begin{equation}\label{sign-distribution}
\pp_{A,\eta}(\sigma)=\prod_{(u,v)\in A}\left(
\frac{\sigma_{(u,v)}+1}{2}-\sigma_{(u,v)}\eta\right),
\end{equation}
and vol stands for the Lebesgue measure on the corresponding euclidean spaces.

\bs \subsection{Some graph--theoretical notations and 
definitions}~\label{subsection-graph-theoretical}\

\ms A path in $\G=(A,V)$ is a sequence $u_0,u_1,\ldots,u_k$, with $u_i\in V$ and 
$(u_i,u_{i+1})\in A$ for each $i$. The length of a path is the number of arrows it 
contains. A cycle is a path $u_0,u_1,\ldots,u_k$, where $u_0=u_k$. We say that the 
vertices $u,v\in V$ are connected if there exists a cycle $u_0,u_1,\ldots,u_k$ such that
$u, v\in \{u_0,u_1,\ldots, u_k\}$. The connected components of $\G$ are the maximal subgraphs 
of $\G$ such that all their vertices are connected. The distance between two vertices 
$u,\ v\in V$, which we denote $\dist(u,v)$, is the length of the shortest directed 
path whose end vertices are $u$ and $v$.

\bs \subsection{The oscillatory subnetwork and the 
asymptotic period}~\label{subsection-the-oscillatory-subnetwork}\

\ms To each couple $(F_{\G,\sigma,T,a},\x)$ we associate a directed graph 
$\G_\osc:=(V_\osc,A_\osc)\subseteq\G$, the oscillatory subnetwork, defined by
\begin{eqnarray*}
A_\osc&:=&\{(u,v)\in A:\ H(\x_u^t-T_{u,v}) \text{ does not converge} \},\\
V_\osc&:=&\{v\in V:\ \exists \ u\in V \text{ such that } 
                             \{(u,v),(v,u)\}\cap A_\osc\neq \emptyset\}.
\end{eqnarray*}

\ms By definition, the oscillatory subnetwork is spanned by all the arrows whose 
activation state $H(\sigma_{u,v}(\x_u-T_{u,v}))$ changes infinitely often. 
This subnetwork is the equivalent, in discrete--time regulatory networks, 
to the dynamical islands introduced in~\cite{Gomez&al06} for 
continuous--time regulatory networks. They are also related to
the clusters of relevant elements in Boolean networks, as they were defined 
in~\cite{Bastolla&Parisi98}. Unlike the differential equations
and the finite state modeling strategies 
of regulatory dynamics~\cite{Glass&Pasternack78, Thieffry07, Thomas73},
our models admit 
oscillations in all network topologies, regardless the distribution of activations and inhibitions through
the graph. This is possible because of the discreteness of time, and continuity of the state variable.
Therefore, oscillatory subnetworks may occur in any discrete--time regulatory network.

\ms The oscillatory subnetwork can be seen as a random variable 
\[
\G_\osc:\F_{a,V}\times [0,1]^V\to \GG_{\subseteq V}:=\{(\bV, \bA): \bV\subset V 
\text{ and } \bA\subset \bV\times \bV\},
\]
from which we can derive other statistical indicators, as for instance its size
$\#V_\osc$, the number ${\rm nc}(\G_\osc)$ of its connected components, and its 
the degree distribution
\[
p_{\G_\osc}(k):=\frac{\#\{v\in V_\osc:\ \Id_\osc(v)=k\}}{\#V_\osc},
\]
where $\Id_\osc(v):=\{u\in V_\osc: (u,v)\in A_\osc\}$.

\ms Another statistical indicator we will consider is the asymptotic period of
a given orbit. Denote by $\per_\tau(F_{\G,\sigma,T,a})$ the set of all 
$F_{\G,\sigma,T,a}$--periodic points of minimal period $\tau$, i.~e., 
\begin{equation}\label{periodic--points}
\per_\tau(F_{\G,\sigma,T,a})=\left\{\x\in[0,1]^V:\ F_{\G,\sigma,T,a}^\tau(\x)=\x \ 
\text{ and } \
F_{\G,\sigma,T,a}^t(\x)\neq \x \ \text{ if }\ t < \tau
\right\}.
\end{equation}
We will say that the asymptotic $F_{\G,\sigma,T,a}$--period of $\x\in [0,1]^V$ 
is $\tau$ if there exists $\y\in\per_\tau(F_{\G,\sigma,T,a})$,
such that 
\[
\lim_{t\to\infty}|F_{\G,\sigma,T,a}^t(\x)-F_{\G,\sigma,T,a}^t(\y)|=0.
\]
We will denote this by $P(F_{\G,\sigma,T,a},\x)=\tau$.

\bs \section{Results}~\label{section-results}

\ms \subsection{The asymptotic period.} \

\ms To each finite set $V$, and $a\in[0,1)$, we associate the sample space 
$\F_{a,V}\times [0,1]^V$, with $\F_{a,V}$ as in 
Equation~\eqref{equation-transformations-set}. 
We supply this sample space with the product sigma--algebra, taking for $[0,1]^V$ and 
$[0,1]^A$, the corresponding Borel sigma--algebras.
The next result concerns the asymptotic period 
$(F_{\G,\sigma,T,a},\x)\mapsto P(F_{\G,\sigma,T,a},\x)$.

\ms
\begin{teorema}\label{theorem-asymptotic-period}
Given a finite set $V$, and $a\in [0,1)$, the function $P:\F_{a,V}\times [0,1]^V\to\nn$
which assigns to $(F_{\G,\sigma,T,a},\x)\in \F_{a,V}\times [0,1]^V$ the value of 
the asymptotic $F_{\G,\sigma,T,a}$--period of $\x$, is a measurable function. 
Furthermore, for each $\eta\in [0,1]$,
\[
\pp_{a,\eta}\left\{(F_{\G,\sigma,T,a},\x)\in \F_{a,V}\times [0,1]^V:\ 
P(F_{\G,\sigma,T,a},\x)<\infty\right\}=1.
\]
\end{teorema}

\ms This expected but nontrivial result says that a random orbit will almost certainly 
approach a periodic orbit, otherwise said, an initial condition almost certainly converges
to a periodic attractor. It is because of this result that we can restrict ourselves 
to the study of periodic orbits, even though the parameter set for which
the orbits are infinite may be uncountable (see~\cite{Coutinho&al06}). We could therefore
redefine the oscillatory subnetwork by saying that, after a transitory, their arrows 
change periodically their activation state.
The proof of this result and some related comments are left to 
Subsection~\ref{subsection-asymptotic-period}.

\bs \subsection{Modularity.} \

\ms In this paragraph we will consider probability measures obtained by 
conditioning on a subnetwork, from a given statistical
ensemble of regulatory networks. Given $\pp_{a,\eta}$ on $\F_{a,V}\times [0,1]^V$, 
and a subnetwork  $\bG:=(\bV,\bA)\in \GG_{\subseteq V}$, we define the analogous
probability measure $\pp_{a,\eta,\bG}$ on $\F_{a,\bV}\times [0,1]^{\bV}$,
such that
\begin{equation}\label{equation-analogous-probability}
\pp_{a,\eta,\bG}\left\{(F_{\G,\bsigma,\bT,a},\y): \bT\in\bar{\I},\ \y\in\bar{\J} \right\}
:=\left\{\begin{array}{ll} 
\pp_{\bA,\eta}(\bsigma)\ \vol(\bar{\I})\ \vol(\bar{\J}) & \text{ if } \G=\bG, \\
                                       0             & \text{ otherwise, }
\end{array}\right.
\end{equation}
for each $\bsigma\in\{-1,1\}^{\bA}$, and all rectangles 
$\bar{\I}\subset [0,1]^{\bA}$ and $\bar{\J}\subset [0,1]^{\bV}$. As before,
vol denotes the corresponding Lebesgue measures. This measure can also be seen as
the marginal of $\pp_{a,V}$ obtained by projecting over $\bV\subset V$, then conditioned 
to have underlying network ${\bG}\in\GG_{\bV}$.

\bs 
\begin{teorema}\label{theorem-modularity}
Fix $a\in[0,1)$, $\eta\in [0,1]$, and a digraph $\bG=(\bV,\bA)\in\GG_{\subseteq V}$ with
$\bV\subsetneq V$. If the digraph distribution $\pp_G$ is such that 
$\pp_G(\G)>0$ for all $\G\in\GG_V$, and if $\pp_{a,\eta}(\G_\osc=\bG) > 0$, then 
we can associated to $\bG$: 
\begin{itemize} 
\item[{\it a)}] a digraph extension $\bG_\ext:=(V,A)\in\GG_V$,
\item[{\it b)}] rectangles 
$\I:=\bI\times\I'\subset [0,1]^{\bA}\times [0,1]^{A\setminus\bA}$ 
and $\J:=\bJ\times\J'\subset [0,1]^{\bV}\times [0,1]^{V\setminus\bV}$, 
\item[{\it c)}] and for each $\bsigma\in \{-1,1\}^{\bA}$, an extension 
$\bsigma_\ext\in \{-1,1\}^A$ such that $\bsigma_\ext|_{\bA}=\bsigma$. 
\end{itemize}
These associated objects satisfy  
\[\G_\osc(F_{\bG_\ext,\bsigma_\ext,T,a},\x) \subseteq\bG\ \ 
\forall\ T\in \I,\ \x\in\J, \text{ and } \bsigma \in \{-1,1\}^{\bA}.\]
Furthermore, $\G_\ext$, $\I$, and $\J$, define surjective transformations
\begin{itemize}
\item[{\it d)}]
$\Phi_A:\bI\times\I'\to [0,1]^{\bA}$, such that 
$\Phi_A(\bT\times T')=\D_{\bA}\bT+\C_{\bA}$ with $\D_{\bA}$ 
diagonal and $\C_{\bA}$ constant, 
\item[{\it e)}] and similarly, 
$\Phi_V:\bJ\times\J'\to [0,1]^{\bV}$ such that 
$\Phi_V(\bx\times \x')=\D_{\bV}\bx+\C_{\bV}$,
with $\D_{\bV}$ diagonal and $\C_{\bV}$ constant. 
\end{itemize}
These transformations satisfy
\[\Phi_V\circ F_{\G,\sigma,T,a}(\x)=F_{\bG,\bsigma,\Phi_A(T),a}\circ\Phi_V(\x),\] 
for all $T\in \I$, and $\x\in\J$.
\end{teorema}

\ms The proof of this result follows from a construction, and it is postponed to 
Section~\ref{section-proofs}. Let us from now expose its interpretation and some of
its consequences. We can divide this Theorem into two parts, the first part 
concerns the stability of the oscillatory subnetworks. 
It says that if an oscillatory subnetwork $\bG:=(\bV,\bA)$ has positive probability to occur, then 
to each signs matrix $\bsigma\in\{-1,1\}^{\bA}$ it corresponds an oscillatory subnetwork 
which is subgraph of  $\bG$. Furthermore, the set of initial conditions and thresholds for which 
a given subgraph of $\bG$ is the oscillatory subnetwork, is a rectangle with nonempty interior. 
Therefore, each of those subgraphs is stable under small changes in the thresholds $T\in [0,1]^A$, 
and in the initial condition $\x\in[0,1]^V$. The size of the maximal perturbation
depends on the position of $(T,\x)$ with respect to the borders of the above mentioned rectangle. 

\ms More interestingly, the second part of the theorem states that if an oscillatory 
subnetwork has positive probability to occur, then the dynamics one can observe 
when restricted to that subnetwork is equivalent to the dynamics supported by the 
subnetwork considered as an isolated system. The equivalence is achieved through the
transformation $\Phi_A$ and $\Phi_V$. It is because of this equivalence that we can talk 
about modularity. Indeed, if we call module any oscillatory subnetwork occurring with 
positive probability, then this theorem can be rephrased by saying that each orbit 
admissible in a module considered as an isolated system, can be achieved, up to a 
change of variables, as the restriction of an orbit in the original system. 
The network extension $\bG_\ext$, the values of the external thresholds 
$\I'\in [0,1]^{\bA\setminus A}$, and the external initial conditions 
$\J'\in [0,1]^{\bV\setminus V}$, can be though as the analogous in our systems,
to the functionality context defined in~\cite{Naldi&al07} for logical networks. 
Using their nomenclature, our theorem ensures that
if an oscillatory subnetwork has a positive probability 
to occur, then there is positive measure set of contexts making this subnetwork
functional.

\ms An interesting consequence of the previous theorem is the following.

\ms \begin{corolario}~\label{corollary-product-periods}
Fix $a\in[0,1)$, $\eta\in [0,1]$, and $\bG\equiv\G_1\cup\G_2\in\GG_{\subseteq V}$, 
with $\G_1:=(V_1,A_1)$ and $\G_2:=(V_2,A_2)$ vertex disjoint and such that $V_1\cup V_2\subsetneq V$. 
If the digraph distribution $\pp_G$ is such that 
$\pp_G(\G)>0$ for all $G\in\GG_V$, and if $\pp_{a,\eta}(\G_\osc\subseteq\bG) > 0$, then 
\[
\pp_{a,\eta}(P=\tau)\geq C \ \sum_{\lcm(\tau_1,\tau_2)=\tau}\pp_{a,\eta,\G_1}(P=\tau_1)\ 
\pp_{a,\eta,\G_2}(P=\tau_2),
\]
with $C>0$ a constant depending on $\pp_{a,\eta}$ and $\bG$.
\end{corolario}

\ms The probabilities in the sum at right hand side of the inequality are analogous probability measures of the kind defined by
Equation~\eqref{equation-analogous-probability}. They are obtained as marginal 
probability corresponding to projections over vertex the sets $V_1$ and $V_2$, then conditioned to have underlying network $\G_1$
and $\G_2$ respectively. 

\ms In~\cite{Bastolla&Parisi98} the authors study the relationship between the modular structure 
and the periods distribution in Boolean networks. The previous corollary establishes such a relation in the framework of discrete--time 
regulatory networks. Our result relates the distribution of the asymptotic period of a network $\G$, to the distribution of the asymptotic 
period of its oscillatory subnetworks considered as isolated systems. 
It establishes in particular that the observable asymptotic periods of $\G$ are the least common multiples of the observable asymptotic 
periods of its oscillatory subnetworks.

\ms \begin{proof}
First notice that, for each $\bG\in \GG_{\subseteq V}$ we have
\[
\pp_{a,\eta}(P=\tau)\geq  
\pp_{a,\eta}(P=\tau|\G_\osc\subseteq\bG)\times \pp_{a,\eta}(\G_\osc\subseteq\bG).
\]
Let $\bG=(\bV,\bA)=(V_1\cup V_2, A_1\cup A_2)$, with $\G_1:=(V_1,A_1)$ and 
$\G_2:=(V_2,A_2)$ vertex disjoint. 
Theorem~\ref{theorem-modularity} ensures the existence of the following objects 
associate to $\bG$:
a) the extension of $\bG$, a digraph $\bG_\ext:=(V,A)\in \GG_V$ such that 
$\bG\subseteq \bG_\ext$, b) two affine surjective transformations 
\[
\Phi_A:\I\to [0,1]^{A_1}\times [0,1]^{A_2} \text{ and } 
\Phi_V:\J\to[0,1]^{V_1}\times [0,1]^{V_2},
\] 
and c) for each $\bsigma\equiv \bsigma^{(1)}\times\bsigma^{(2)}
\in\{-1,1\}^{A_1}\times \{-1,1\}^{A_2}$,
an extension $\bsigma_\ext\in\{-1,1\}^{A}$ such that $\sigma|_{A_1\cup A_2}=\bsigma$. 
From the second part of the theorem it follows that 
\begin{equation}~\label{equation-factorization-mappings}
        \Phi_V\left(F^{t}_{\bG_\ext,\bsigma_\ext,T,a}(\x)\right) 
         = F^{t}_{\G_1,\bsigma^{(1)},\bT^{(1)},a}\left(\y^{(1)}\right)
    \times F^{t}_{\G_2,\bsigma^{(2)},\bT^{(2)},a}\left(\y^{(2)}\right),
\end{equation}
for all $t\in\nn$ and $(T,\x)\in \I\times\J$. Here we have used 
$\Phi_V(T)\equiv\bT^{(1)}\times\bT^{(2)}\in [0,1]^{A_1}\times [0,1]^{A_2}$ and  
$\Phi_V(\x)\equiv \y^{(1)}\times\y^{(2)}\in [0,1]^{V_1}\times [0,1]^{V_2}$.
For $i=1, 2$, and each $\bsigma^{(i)}\in \{-1,1\}^{A_i}$ and $\tau_i \in \nn$, 
let us define
\[
P_{\G_i,\bsigma^{(i)},\tau_i}:=\left\{\left(\bT^{(i)},\y^{(i)}\right)
\in [0,1]^{A_i}\times [0,1]^{V_i}:
\ P\left(F_{\G_i,\bsigma^{(i)}, \bT^{(i)},a},\y^{(i)}\right)=\tau_i\right\}.
\]
Defining $\Phi:\I\times\J\to [0,1]^{A}\times[0,1]^V$ such that 
$\Phi(T,\x):=\Phi_A(T)\times \Phi_V(\x)$, it follows 
from~\eqref{equation-factorization-mappings} 
that
\begin{eqnarray*}
\pp_{a,\eta}\left(P=\tau|\G_\osc \subseteq \bG\right)
&\geq & \frac{\pp_G(\bG_\ext)}{\pp_{a,\eta}(\G_\osc \subseteq\bG)}
\sum_{\bsigma=\bsigma^{(1)}\times \bsigma^{(2)}}\pp_{A,\eta}(\bsigma_\ext)\\
&     &\hskip 40pt \times \sum_{\lcm(\tau_1,\tau_2)=\tau}
\vol\circ\Phi^{-1}\left(P_{\G_1,\bsigma^{(1)},\tau_1}\times P_{\G_2,\bsigma^{(2)},\tau_2}\right),
\end{eqnarray*}
where vol denotes Lebesgue measure in $[0,1]^A\times[0,1]^V$. 
Let us now use the decomposition 
$\I=\bI\times\I'\subset [0,1]^{\bA}\times [0,1]^{A\setminus\bA}$ 
and $\J=\bJ\times\J'\subset [0,1]^{\bV}\times [0,1]^{V\setminus\bV}$. 
According to Theorem~\ref{theorem-modularity}, we have
\[
\Phi((\bT\times T')\times(\bx\times\x'))= \D(\bT\times\bx)+\C
\]
where $\D$ is a diagonal linear bijection, and $\C$ is a constant. With this, 
\begin{eqnarray*}
\pp_{a,\eta}\left(P=\tau|\G_\osc=\bG\right)
&\geq & \frac{\pp_G(\bG_\ext)\ \vol(\I'\times \J')}{|\D|\ \pp_{a,\eta}(\G_\osc\subseteq\bG)}
\sum_{\bsigma=\bsigma^{(1)}\times \bsigma^{(2)}}\pp_{A,\eta}(\bsigma_\ext)\\
&     &\hskip 40pt \times \sum_{\lcm(\tau_1,\tau_2)=\tau} 
\vol_1\left(P_{\G_1,\bsigma^{(1)},\tau_1}\right)\times 
      \vol_2\left(P_{\G_2,\bsigma^{(2)},\tau_2}\right),
\end{eqnarray*}
where vol denotes the Lebesgue measure in 
$[0,1]^{A\setminus \bA}\times [0,1]^{V\setminus\bV}$, 
vol${}_i$ denotes the Lebesgue measure in $[0,1]^{A_i}\times [0,1]^{V_i}$ for $i=1,2$, 
and $|\D|$ denotes the determinant of the transformation $\D$. 
From here, and taking into account the definition in 
Equation~\eqref{equation-analogous-probability}, we have
\begin{eqnarray*}
\pp_{a,\eta}\left(P=\tau|\G_\osc=\bG\right)
&\geq & \frac{\pp_G(\bG_\ext)\ \vol(\I'\times \J')
               \min(\eta,1-\eta)^{\#(A\setminus\bA)}}{|\D|\ \pp_{a,\eta}(\G_\osc\subseteq\bG)} \\
&    &\hskip 30pt \times \sum_{\bsigma=\bsigma^{(1)}\times \bsigma^{(2)}}       
\pp_{A_1,\eta}\left(\bsigma^{(1)}\right)\times\pp_{A_2,\eta}\left(\bsigma^{(2)}\right)\\
&     &\hskip 50pt \times \sum_{\lcm(\tau_1,\tau_2)=\tau} 
\vol_1\left(P_{\G_1,\bsigma^{(1)},\tau_1}\right)\times 
 \vol_2\left(P_{\G_2,\bsigma^{(2)},\tau_2}\right)\\
&\geq & \frac{\pp_G(\bG_\ext)\ \vol(\I'\times \J')\ 
\min(\eta,1-\eta)^{\#(A\setminus\bA)}}{|\D|\ \pp_{a,\eta}(\G_\osc\subseteq\bG)}\\
&     &\hskip 40pt \times \sum_{\lcm(\tau_1,\tau_2)=\tau} 
\pp_{a,\eta,\G_1}(P=\tau_1)\ \pp_{a,\eta,\G_2}(P=\tau_2),
\end{eqnarray*}
and the result follows with 
$C:=\pp_G(\bG_\ext)\times \vol(\I'\times \J')\times 
       \min(\eta,1-\eta)^{\#(A\setminus\bA)}/|\D|$.
\end{proof}

\bs \subsection{Sign Symmetry.} \

\ms The next result illustrates how a structural constrain in the underlying random 
network manifests in the distribution of the statistical indicators. Here, the fact that 
the underlying random digraph does not admit cycles of odd length, implies a symmetry on 
the statistical indicators with respect to the change of sign in the interactions. 
Let us remind that for each $\eta\in[0,1]$ and $a\in[0,1)$, the probability measure 
$\pp_{a,\eta}$ on $\F_{a,V}\times[0,1]^V$ which defines a statistical ensemble, is obtained 
from a fixed distribution $\pp_G$ over the finite set of all directed graphs with 
vertices in $V$. We have the following.

\ms \begin{proposicion}~\label{proposition-symmetry} 
If the distribution $\pp_G$ over the set $\GG_V$ of all directed graphs with 
vertices in $V$ is such that 
$\pp_G\{\G\in \GG_V: \G \text{ admits a cycle of odd length }\}=0$,
then
\begin{eqnarray*}
\pp_{a,\eta}\left(P=\tau\right)&=&
                   \pp_{a,1-\eta}\left(P=\tau\right), \ \ \forall \tau\in\nn,\\
\pp_{a,\eta}\left(\G_\osc=\bG\right)&=&
   \pp_{a,1-\eta}\left(\G_\osc=\bG\right)\ \ \forall \ \bG\in\GG_{\subseteq V}, 
\end{eqnarray*}
for each $\eta\in[0,1]$.
\end{proposicion}

\ms Statistical ensembles of regulatory networks whose underlying digraphs are 
random trees (as defined in Subsection~\ref{subsection-random-digraphs} below) 
clearly satisfy the hypothesis of Proposition~\ref{proposition-symmetry}. 
Random subnetworks of the square lattice also have this property.
Therefore, for regulatory dynamics over those kind of digraphs, the statistical 
indicators we consider are left invariant under the symmetry $\eta\mapsto 1-\eta$. 

\ms In the proof of Proposition~\ref{proposition-symmetry}, we will need the following.

\ms 
\begin{lema}~\label{lemma-notin-delta} 
For each $\G:=(V,A)\in\GG_V$, and $\sigma\in\{-1,1\}^V$ we have,
\[
\vol\left\{(T,\x)
\in[0,1]^A\times [0,1]^V:\ \left\{F^t_{\G,\sigma,T,a}(\x):\ t\in \nn\right\}\cap 
\Delta_T=\emptyset\right\}=1.
\] 
where $\Delta_T:=\left\{\x\in[0,1]^V:\ \x_u=T_{u,v}\ \text{for some} (u,v)\in A\right\}$ is the discontinuity 
set of the piecewise constant part of $F_{\G,\sigma,T,a}$.
\end{lema}

\ms This lemma directly follows from the arguments developed below, in  
paragraphs~\ref{subsubsection-density-argument} and~\ref{subsubsection-complexity-argument},
inside the proof of Theorem~\ref{theorem-asymptotic-period}.

\ms {\it Proof of Proposition~\ref{proposition-symmetry}}.
First notice that, for each $\G=(A,V)\in \GG_V$ and $\sigma\in\{-1,1\}^A$ 
we have 
\begin{eqnarray*}
\pp_{A,1-\eta}(-\sigma)&=&
\prod_{(u,v)\in A}\left(
          \frac{-\sigma_{(u,v)}+1}{2}+\sigma_{(u,v)}(1-\eta)\right)\\
                      &=&
\prod_{(u,v)\in A}\left(
          \frac{\sigma_{(u,v)}+1}{2}-\sigma_{(u,v)}\eta\right)=\pp_{A,\eta}(\sigma). 
\end{eqnarray*}
Let us suppose that all the cycles in $\G$ have even length, and
for each connected component of $\bG:=(\bA,\bV)\subseteq\G=(A,V)\in \GG_V$,
choose a pivot vertex $\bu\in \bV$. With this define $\Psi_A:[0,1]^A\to [0,1]^A$ and 
$\Psi_V:[0,1]^V\to[0,1]^V$ such that
\[
\Psi_A(T)_{(u,v)}:=\left\{\begin{array}{ll}
                    1-T_{(u,v)} & \text{ if } \dist(u,\bu) \in 2\nn\\
                    T_{(u,v)}   & \text{ otherwise }
                    \end{array}\right. \hskip 5pt
\Psi_V(\x)_u:=\left\{\begin{array}{ll}
                    1-\x_u & \text{ if } \dist(u,\bu)\in 2\nn\\
                    \x_u   & \text{ otherwise, }
                    \end{array}\right.
\]    
for each $u,v\in \bV$. Here $d_\G$ denotes vertex distance induced by the digraph
$\G$, as defined in Subsection~\ref{subsection-graph-theoretical}.               
Both $\Psi_A$ and $\Psi_V$ are affine isometries, therefore they preserve 
the Lebesgue measure in $[0,1]^A$ and $[0,1]^V$ respectively. 
If $\dist(v,\bu)$ is even, then $\Psi_V(\x)_v=1-\x_v$, and for each
$u$ such that $(u,v)\in A$, $\Psi_V(\x)_u=\x_u$ and
$\Psi_A(T)_{(u,v)}=T_{(u,v)}$. In this case we have
\begin{eqnarray*}
F_{\G,-\sigma,\Psi_A(T),a}(\Psi_V(\x))_v&=&a(1-\x_v)+
        \frac{1-a}{\Id(v)}\sum_{u\in V: (u,v)\in A}
     H\left(-\sigma_{(u,v)}\left(\x_u-T_{(u,v)}\right)\right)\\
                    &=&a(1-\x_v)+
     \frac{1-a}{\Id(v)}\sum_{u\in V: (u,v)\in A} \left(1 -
     H\left(\sigma_{(u,v)}\left(\x_u-T_{(u,v)}\right)\right)\right)\\
                    &=&1-F_{\G,\sigma,T,a}(\x)_v=
                                   \Psi_V(F_{\G,\sigma,T,a}(\x))_v,
\end{eqnarray*} 
for each $\x\not\in\Delta_T$. On the other hand, if $\dist(v,v_{\rm pivot})$ is odd, 
we have $\Psi_V(\x)_v=\x_v$, and for each $u$ such that $(u,v)\in A$, 
$\Psi_V(\x)_u=1-\x_u$ and $\Psi_A(T)_{(u,v)}=1-T_{(u,v)}$. 
In this case 
\begin{eqnarray*}
F_{\G,-\sigma,\Psi_A(T),a}(\Psi_V(\x))_v&=&
a\x_v +\frac{1-a}{\Id(v)}\sum_{u\in V: (u,v)\in A}
H\left(-\sigma_{(u,v)}\left((1-\x_u)-\left(1-T_{(u,v)}\right)\right)\right)\\
               &=&a\x_v+\frac{1-a}{\Id(v)} \sum_{u\in V: (u,v)\in A}
               H\left(\sigma_{(u,v)}\left(\x_u-T_{(u,v)}\right)\right) \\
              &=&F_{\G,\sigma,T,a}(\x)_v=\Psi_V(F_{\G,\sigma,T,a}(\x))_v,
\end{eqnarray*} 
for each $\x\not\in\Delta_T$. Taking into account 
Lemma~\ref{lemma-notin-delta}, $F^t_{\G,-\sigma,\Psi_A(T),a}(\Psi_V(\x))=\Psi_V(F^t_{\G,\sigma,T,a}(\x))$
for all $t\in\nn$, for almost all $(\x,T)\in [0,1]^A\times[0,1]^V$. 
Since $\Phi_V$ is an isometry, we immediately have 
\[
\vol\left\{(T,\x)\in[0,1]^A\times[0,1]^A:\ 
P(F_{\G,-\sigma,\Psi_A(T),a},\Psi_V(\x))=P(F_{\G,\sigma,T,a},\x)\right\}=1,
\]
for all $\G:=(V,A)\in\GG_V$ with no odd cycles, and every $\sigma\in\{-1,1\}^A$.  
Since both $\Psi_A$ and $\Psi_V$ preserve the respective Lebesgue measure, and 
$\pp_G( \GG \text{ admits an odd cycle})=0$, we obtain
\begin{eqnarray*}
\pp_{a,\eta}(P=\tau)&=&\sum_{\G\in\GG_V}\pp_G(\G)\sum_{\sigma\in\{-1,1\}^A}
\pp_{A,1-\eta}(-\sigma)\\ 
& &\hskip 105pt\times \vol\left\{(T,\x):\ 
P(F_{\G,-\sigma,\Psi_A(T),a},\Psi_V(\x))=\tau\right\}\\
&=&\sum_{\G\in\GG_V}\pp_G(\G)\sum_{\sigma\in\{-1,1\}^A}
\pp_{A,1-\eta}(-\sigma)\ 
\vol\left\{(T,\x):\ 
P(F_{\G,-\sigma,T,a},\x)=\tau\right\}\\
&=&\pp_{a,1-\eta}(P=\tau),
\end{eqnarray*}
for all $\tau\in\nn$, and this concludes the proof of the first claim in the proposition.

\ms For the second claim, if $F^t_{\G,\sigma,T,a}(\x)\notin\Delta_T$, then
\begin{eqnarray*}
 H\left(\sigma_{(u,v)}(F^t_{\G,\sigma,T,a}(\x)_u-T_{(u,v)})\right)&=&
H\left(-\sigma_{u,v}\left(\Psi_V(F^t_{\G,\sigma,T,a}(\x))_u-
                                           \Psi_A(T)_{(u,v)}\right)\right)\\
                                                                   &=&
H\left(-\sigma_{u,v}\left(F^t_{\G,-\sigma,\Psi_A(T),a}(\Psi_V(\x))_u-
             \Psi_A(T)_{(u,v)}\right)\right).
\end{eqnarray*}
Therefore, by 
Lemma~\ref{lemma-notin-delta},
\[
\vol\left\{(T,\x)\in[0,1]^A\times[0,1]^A:\ 
\G_\osc(F_{\G,-\sigma,\Psi_A(T),a},\Psi_V(\x))=\G_\osc(F_{\G,\sigma,T,a},\x)\right\}=1,
\]
for all $\G:=(V,A)\in\GG_V$ all of whose cycles have even length, and every  
$\sigma\in\{-1,1\}^A$. From here, taking into account that 
$\pp_G(\G \text{ admits an odd cycle})=0$, and the fact 
that that $\Psi_A$ and $\Psi_V$ preserve the Lebesgue measure, we obtain
\begin{eqnarray*}
\pp_{a,\eta}(\G_\osc=\bG)&=&\sum_{\G\in\GG_V}\pp_G(\G)\sum_{\sigma\in\{-1,1\}^A}
\pp_{A,1-\eta}(-\sigma)\\ 
&  &\hskip 105pt \times \vol\left\{(T,\x):\ 
\G_\osc(F_{\G,-\sigma,\Psi_A(T),a},\Psi_V(\x))=\bG\right\}\\
&=&\sum_{\G\in\GG_V}\pp_G(\G)\sum_{\sigma\in\{-1,1\}^A}
\pp_{A,1-\eta}(-\sigma)\ 
\vol\left\{(T,\x):\ 
\G_\osc(F_{\G,-\sigma,T,a},\x)=\bG\right\}\\
&=&\pp_{a,1-\eta}(\G_\osc=\bG),
\end{eqnarray*}
for all $\bG\in\GG_{\subseteq V}$, and the proof is completed.

\endproof

\bs \section{Examples}\label{section-examples}

\ms In this section we present, as a matter of illustration, two families of random digraphs 
which we have explored numerically. A detailed numerical study, which requires massive 
calculations, is out of the purpose of the present work, and is left for a future research. 
Instead, in this section we make some general observations suggested by a preliminary numerical 
exploration of these examples.

ù\bs\subsection{Two Families of Random digraphs}~\label{subsection-random-digraphs}\

\ms We consider two kind of distributions $\pp_G$ on the set 
$\GG_V:=\{(V, A): A\subset V\times V\}$ of all directed graph with vertex set $V$. 
On one hand we have a directed version of the classical Erdo\"s--R\'enyi 
ensemble, which we define as follows.

\ms Fix a vertex set $V$ and $p\in (0,1)$, and consider the probability 
distribution $\pp_{p}:\GG_{V}\to [0,1]$ such that $\pp_{p}(A)=p^{\#A}$. 
According to this, a directed graph contains the arrow $(u,v)$ 
with probability $p$, and the inclusion of different vertices are 
independent and identically distributed random variables. 
The typical Erdo\"s--R\'enyi graph is statistically homogeneous, 
thus suited for a mean field treatment. Most of the rigorous results 
concerning random graphs refer to these kind of models (see~\cite{Bollobas01, Durrett07} 
and references therein).

\ms The other family of random digraphs we will refer to derives from the famous 
Barab\'asi--Albert model of scale--free random graphs, whose popularity relies on  
the ubiquity of the scale--free property~\cite{Barabasi&Oltvai04, Jeong&al00}. 
A dynamical construction of scale--free graphs was proposed by Barab\'asi and Albert 
in~\cite{Barabasi&Albert99} (see~\cite{Durrett07} for a more rigorous presentation). 
Their model incorporates two key ingredients: continuous growth and 
preferential attachment. We implemented their construction as follows.
For $n=0$, we take $\G_0=K_{m_0}$, the complete simple undirected graph in $m_0$ 
vertices. Then, for each $n\geq 1$, a new vertex $v_{n+1}$ is added to the graph 
$\G_n:=(E_n,V_n)$.
This new vertex form new edges with randomly chosen vertices in $V_n$. The 
probability for $v\in V_n$ to be chosen is proportional to its degree in $\G_0$, i.~e, 
\[
\{v_{n+1},v\}\in E_n \text{ with probability } 
p_n(v):=\frac{\#\{ u\in V_n: \ \{u,v\}\in E_n\}}{\#E_n}.
\] 
At the $(n+1)$--th step we obtain a graph $\G_{n+1}$ with an increased vertex set 
$V_{n+1}:=V_n\cup \{v_{n+1}\}$ and an enlarged edge set $E_{n+1}$. 
This random iteration continues until a predetermined number $N$ of vertices is 
obtained. In this scheme we are able to add more than one edge at each iteration. 
If on the contrary, we allow only one new edge at each time step, and we 
start with $m_0=2$ vertices choosing at the $1$--th step one of the two preexisting 
vertices to form a new edge with with probability 1/2, then the resulting 
graph would be a random scale--free tree~\cite{Durrett07}. In any case,
the final graph $\G_{N-m_0}$ is turned into a directed graph $\G=(V,A)$ 
by randomly assigning a direction to each edge $\{u,v\}\in E_{N-m_0}$. In our case 
we choose one of the two possible directions with probability $1/4$,
and give probability $1/2$ to the choice of both directions.  
This procedure generates a probability distribution $\pp_{{\rm BA}}$ in 
$\GG_{V}:=\{(V, A): A\subset V\times V\}$. 

\bs \subsection{Erdo\"s--R\'enyi}~\label{subsection-erdos-numeric}\

\ms We generated random digraphs of the Erdo\"s--R\'enyi kind with 100 vertices,
and probability of connection $p \in \{0.2, 0.4, 0.6, 0.8\}$.
For $p$ fixed and each  
$\eta \in \{0,0.1,0.2,\ldots,1\}$, we build 20 graphs, to which we randomly and
uniformly assign thresholds, and signs with probability $\pp_{A,\eta}$. 
By taking $a \in \{0.0,0.1,0.2,\ldots,0.8\}$, we obtain 7920 regulatory 
networks, 20 for each triplet 
$(p,a,\eta) \in\{0.2,0.4,0.6,0.8\}\times\{0.0,0.1,\ldots,0.8\}\times\{0,0.1,0.2,\ldots,1\}$.
Finally, for each one of those dynamical systems we run 30 initial conditions randomly 
taken from $[0,1]^{100}$, and we let the system evolve in order to determine 
the corresponding oscillatory subnetwork and the asymptotic period.

\ms The characteristics of the oscillatory subnetwork depend strongly on the parameter
$a$, and less noticeably on $\eta$. We have observed that the mean size of the oscillatory 
subnetwork decreases with $a$ in a seemingly monotonous way, and does not appear to depend 
on $\eta$. The mean number of connected components never exceeded 2. 
A crude inspection of the degree distributions of the oscillatory subnetworks suggests 
that they form an ensemble of the Erdo\"s--R\'enyi type, nevertheless, any conclusion in this 
direction requires a more detailed numerical study.
Finally, the mean value of the asymptotic period seems to grow with both $a$ and $\eta$. 
It increases monotonously with $a$, while for each fixed $a$, it reaches a local maximum 
around $\eta=0.5$.

\begin{figure}[h!]
\begin{center}
\begin{tabular}{cc}
\includegraphics[width=0.5\textwidth]{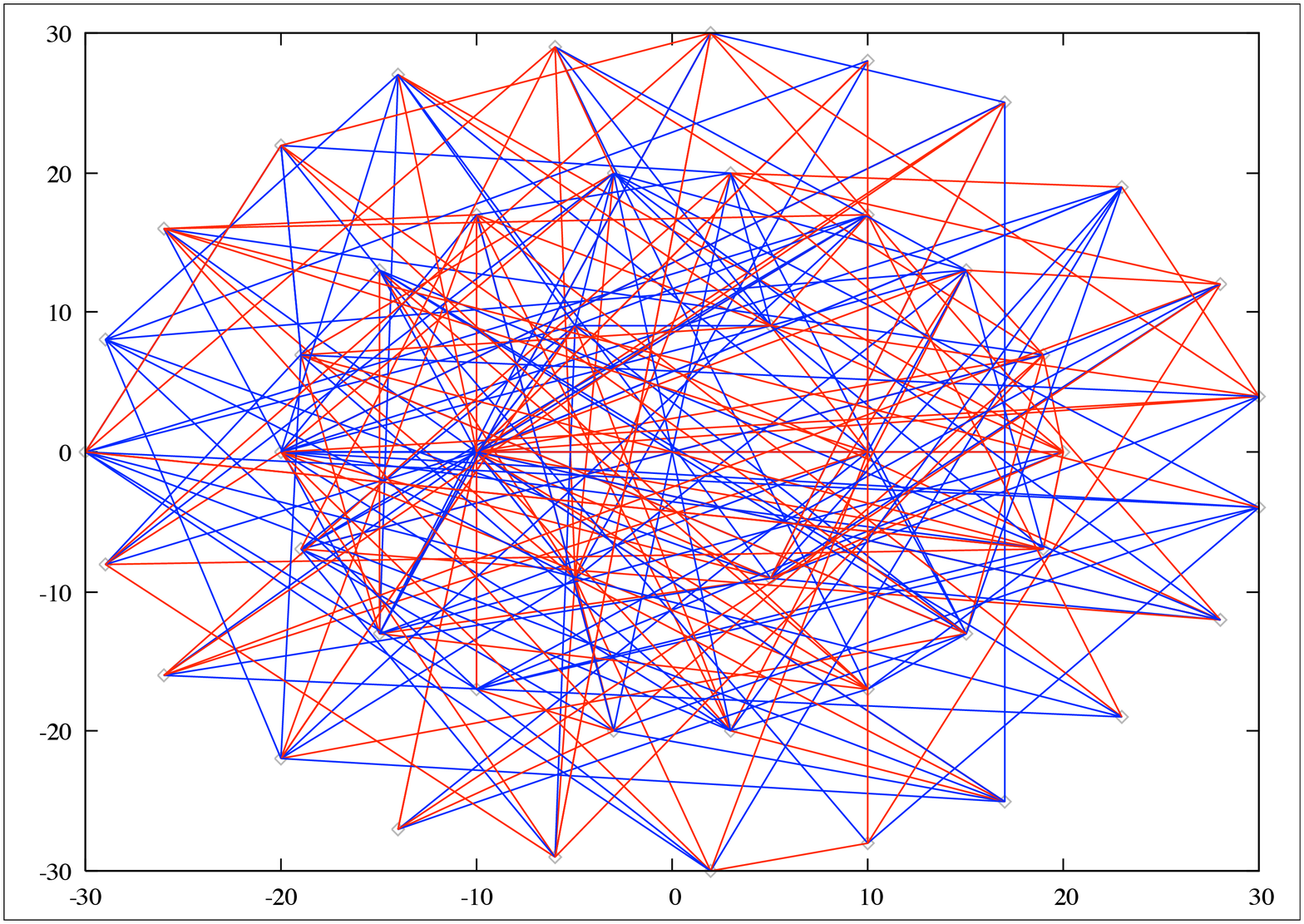} &
\includegraphics[width=0.5\textwidth]{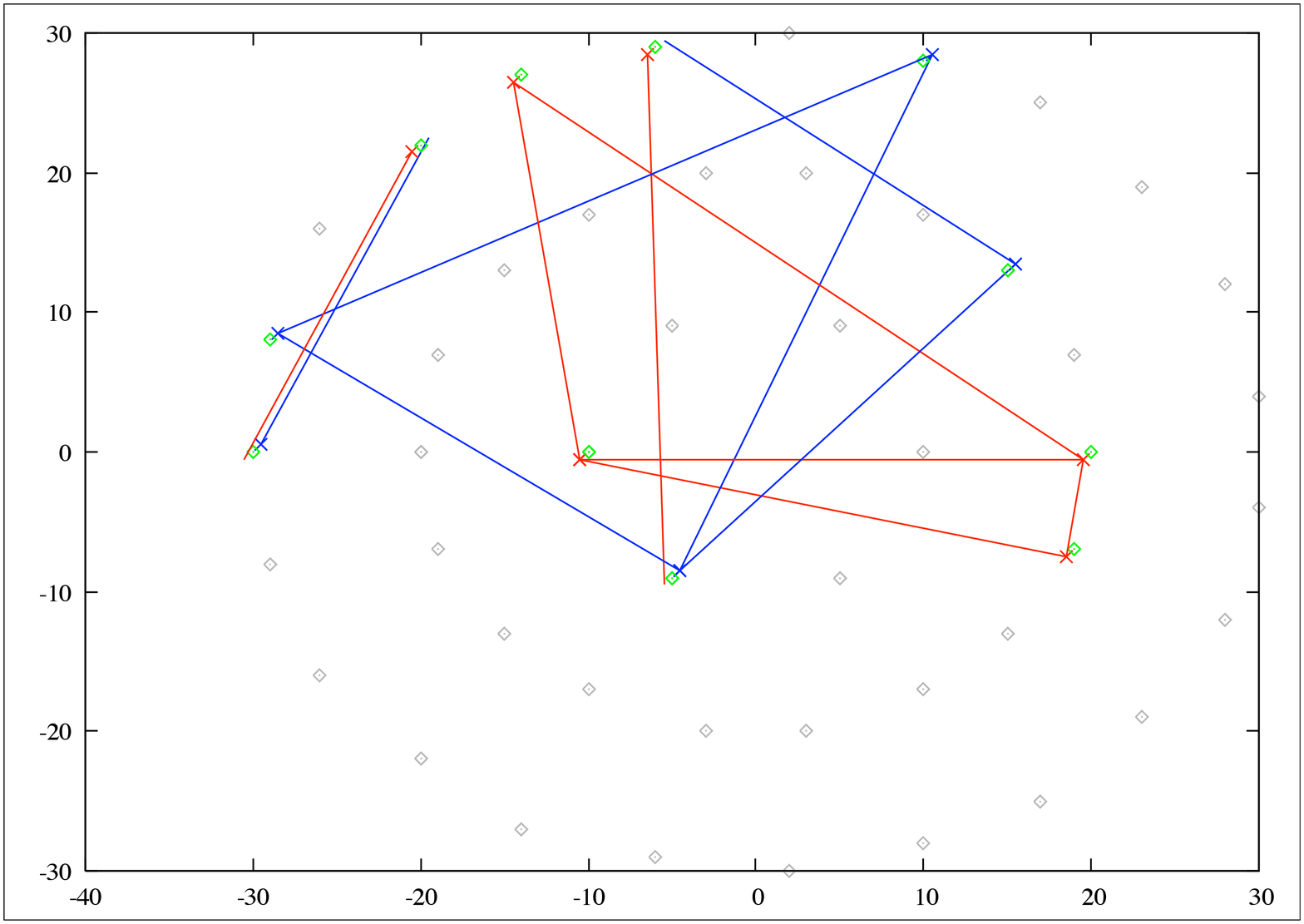}
\end{tabular}
\end{center}
\caption{An Erdo\"s--R\'enyi network and corresponding oscillatory subnetwork. 
The underlying graph was built with a connection probability $p=0.1$, and the 
oscillatory subnetwork corresponds to a random initial condition, 
with random threshold values, and contraction rate $a=0.2$. Signs are codified 
by colors: blue for +1, and red for -1. They appear with equal probability $\eta=0.5$. 
In the oscillatory subnetwork, the direction of an arrow is indicated by marking its
head vertex with a $\times$.}
\label{figure-erdos-graphs}
\end{figure}

\ms As a matter of illustration, in Figure~\ref{figure-erdos-graphs} we show 
two digraph with vertex set $V=\{1,2,\ldots,50\}$.
At the left side we draw a random digraph $\G\in\GG_V$ according to the distribution 
$\pp_{p}$, with $p=0.1$. The signs of the arrows correspond to a random choice according to
$\pp_{A,\eta}$, with $\eta=0.5$, i.~e., both signs have the same probability to be chosen. 
The digraph at right is the oscillatory subnetwork $\bG=\G_\osc(F_{\G,\sigma,T,a},\x)$, 
corresponding to a certain choice of thresholds $T\in [0,1]^A$ and initial 
condition $\x\in[0,1]^V$. The contraction rate was $a=0.2$. In the oscillatory 
subnetwork, the direction of an arrow $(u,v)$ is indicated by marking head $v$ with 
a $\times$ sign. 

\bs \subsection{Barab\'asi--Albert}~\label{subsection-barabasi-numeric}\

\ms Following the procedure described in~\ref{subsection-random-digraphs} we generate,
starting from a complete graph in $m_0=5$ vertices, ensembles of random digraphs of the
Barab\'asi--Albert kind of size $N=100$. For each $\eta \in \{0,0.1,0.2,\ldots,1\}$, 
we generate 20 random graphs, then we turn them into directed graphs by
choosing one of the two possible directions with probability $1/4$,
and both directions with probability $1/2$. Finally we randomly assign thresholds, 
and signs with probability $\pp_{A,\eta}$, to each one of the 20 resulting digraphs. 
By ranging $a \in \{0.0,0.2,0.4,0.6\}$, we obtain 880 regulatory networks, 20 for each 
couple $(a,\eta) \in \{0.0,0.2,0.4,0.6\}\times \{0,0.1,0.2,\ldots,1\}$.
Then, for each one of those dynamical systems, we randomly select 30 initial conditions in
$[0,1]^{100}$, and let the system evolve in order to determine the corresponding 
oscillatory subnetwork and the asymptotic period.

\ms For this family of random networks, the mean size of the oscillatory 
subnetwork behaves in a similar way as for the Erdo\"s-R\'enyi family, i.~e., 
it increases nonmonotonously with $\eta$, with a local maximum at $\eta=0.5$, 
and decreases monotonously with $a$. The mean number of connected components 
(ranging from 0 to 10) shows a similar behavior with respect to both $\eta$ and $a$, with
a very noticeable increase around $\eta=0.5$. On the other hand, the distribution of 
the asymptotic period becomes broader as we approach
$\eta=0.5$, and for each value of $\eta$, the mean value of the asymptotic 
period grows with $a$. It is worth to notice that these distributions become 
broader as $a$ increases, and at the same time the number of connected components 
grows. This behavior is consistent with the proliferation of periods predicted 
by Corollary~\ref{corollary-product-periods}, in the case of a statistical 
ensemble whose oscillatory subnetworks have many connected components.

\begin{figure}[h!]
\begin{center}
\begin{tabular}{cc}
\includegraphics[width=0.5\textwidth]{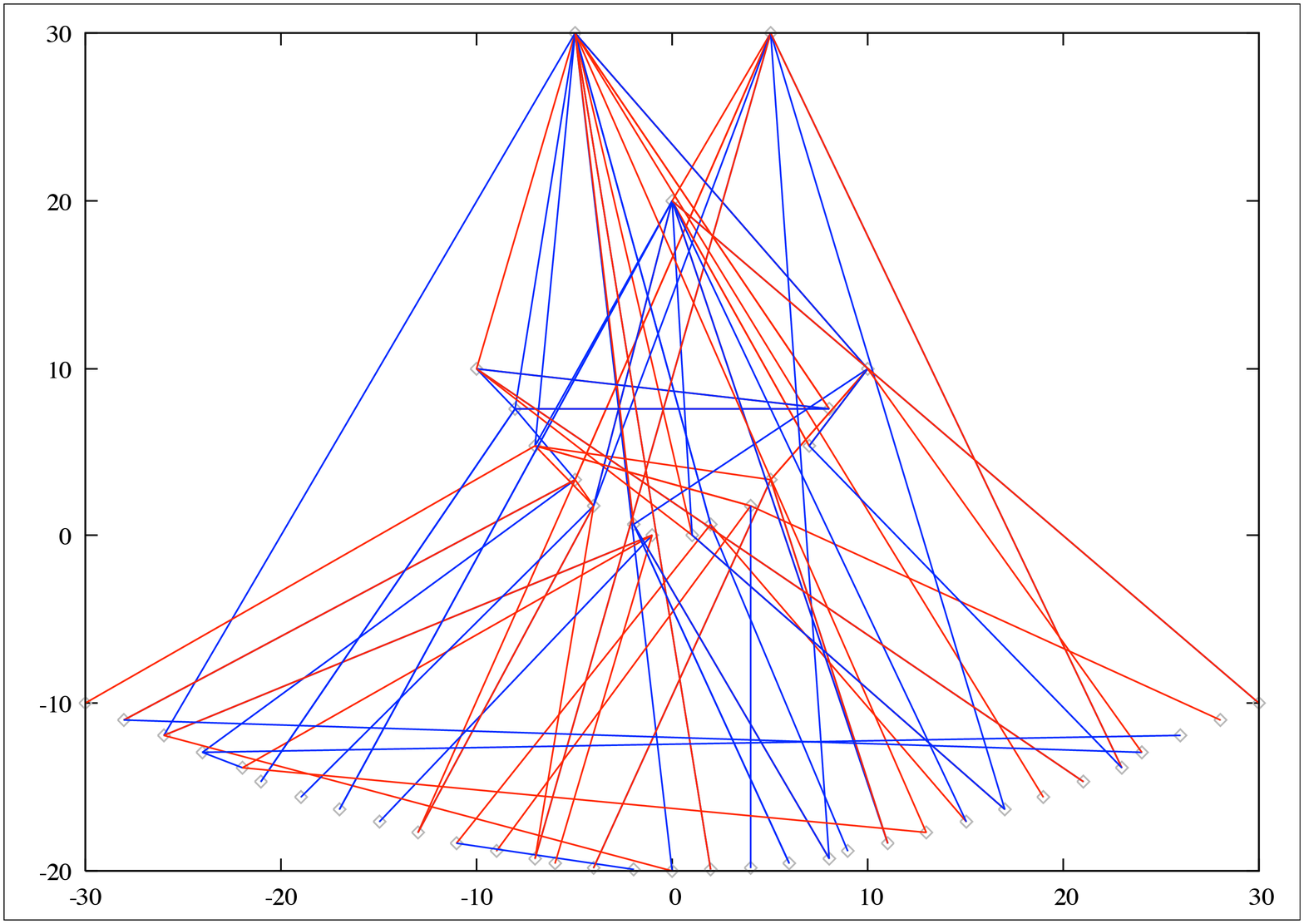} &
\includegraphics[width=0.5\textwidth]{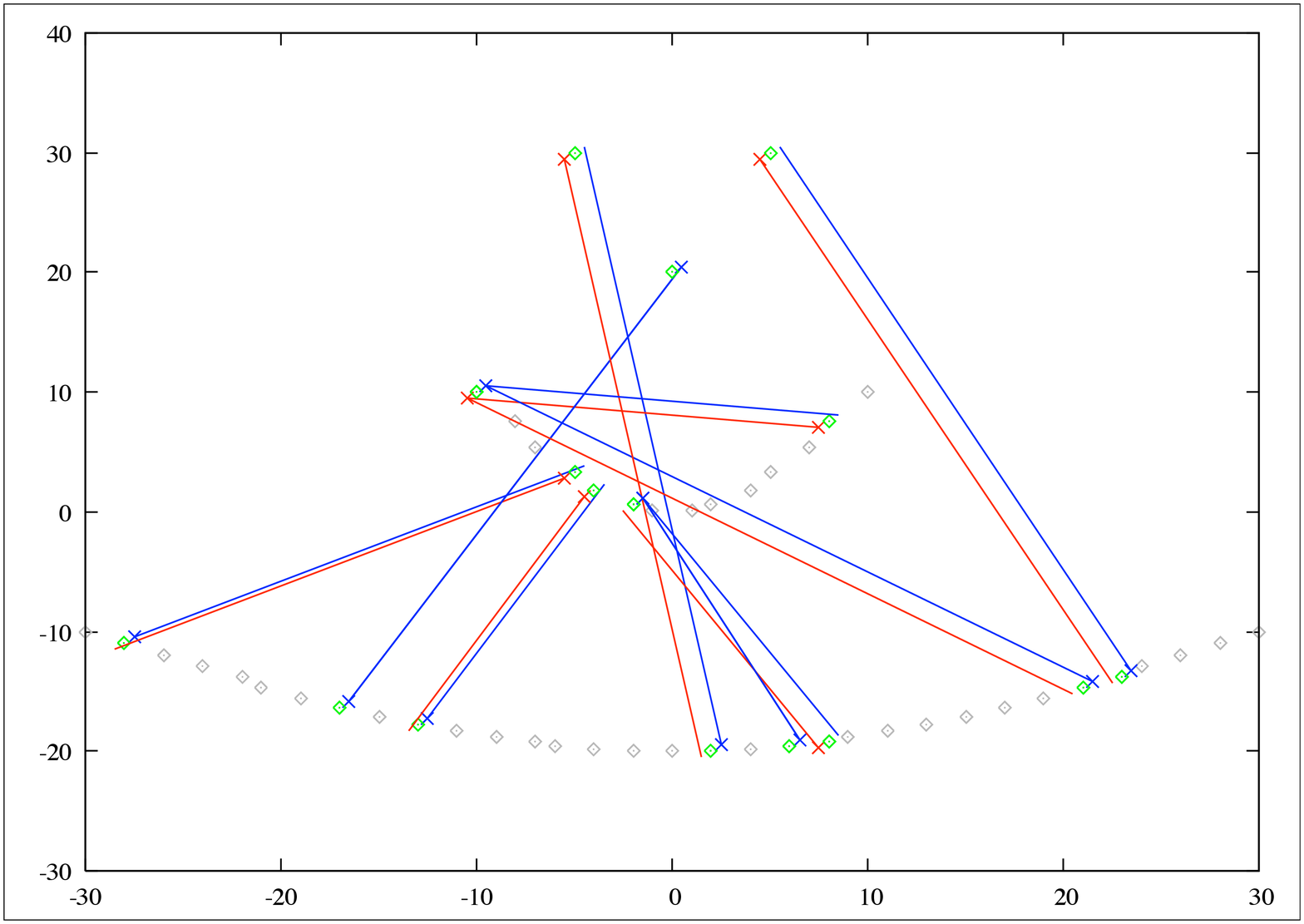}
\end{tabular}
\end{center}
\caption{A Barab\'asi--Albert network and corresponding oscillatory subnetwork. 
The underlying graph was built according to the procedure described in 
Paragraph~\ref{subsection-random-digraphs}, starting with $m_0=3$ vertices. 
The oscillatory subnetwork corresponds to a random initial condition, with random 
threshold values, and contraction rate $a=0.2$. 
Signs are codified by colors: blue for +1, and red for -1. They appear 
with equal probability $\eta=0.5$. 
In the oscillatory subnetwork, the direction of the arrows is indicated by marking 
the head with a $\times$.}
\label{figure-barabasi-graphs}
\end{figure}

\ms In Figure~\ref{figure-barabasi-graphs} we show two digraph with vertex 
set $V=\{1,2,\ldots,50\}$.
At the left side we draw a random digraph $\G\in\GG_V$ according to the distribution 
$\pp_{{\rm BA}}$. The sign of the arrows is randomly chosen by using the distribution 
$\pp_{A,\eta}$ with $\eta=0.5$, i.~e., both signs have the same 
probability. The digraph at right is the oscillatory subnetwork 
$\bG=\G_\osc(F_{\G,\sigma,T,a},\x)$, determined by a random choice of thresholds 
$T\in [0,1]^A$ and initial condition $\x\in[0,1]^V$, with contraction rate $a=0.2$. 
In the oscillatory subnetwork, the direction of the arrows is indicated by marking the 
head with a $\times$. 

\begin{figure}[h!]
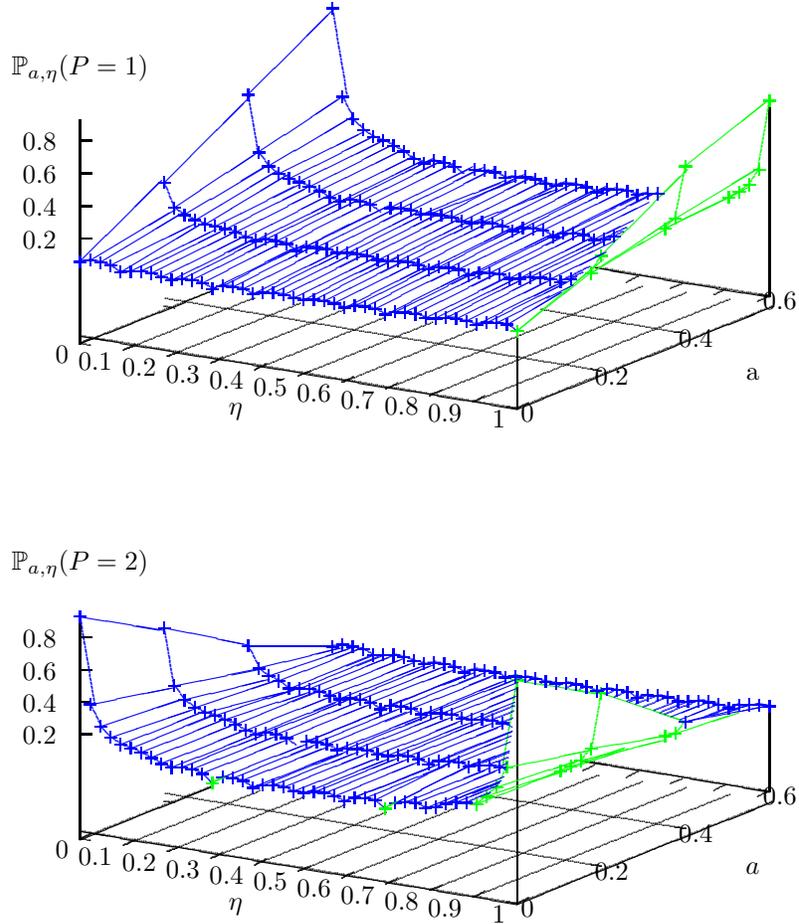

\begin{center}
\begin{tabular}{c}
\scalebox{1}{\input{periodo1B3D.latex}}
\\
\scalebox{1}{\input{periodo2B3D.latex}} \\
\end{tabular}
\end{center}
\caption{At the top, the probability to approach a fixed point 
in the Barab\'asi--Albert ensembles, as a function of $a$ and $\eta$. 
At the bottom, the probability for the asymptotic period to be equals 2, 
also as a function of $a$ and $\eta$, for the same ensembles.}
\label{figure-tonyo}
\end{figure}

\ms In Figure~\ref{figure-tonyo} we show the behavior of the probability to 
approach a fixed point, and the probability for the asymptotic period to be equals 2, 
both as a function of $a$ and $\eta$. This picture summarizes the statistics 
of the 26400 orbits generated on Barab\'asi--Albert like networks, 
and presents the symmetry $\eta\to 1-\eta$ predicted by 
Proposition~\ref{proposition-symmetry}.
This symmetry is expected in an ensemble of regulatory network whose underlying 
digraphs have no odd length cycles. Though the underlying networks of the 
experiments we performed have both odd and even length cycles, there is a larger 
proportion of even lengths amongst the cycles of small length.  

\bs \section{Proofs}~\label{section-proofs} \

\ms \subsection{Proof of 
Theorem~\ref{theorem-asymptotic-period}}~\label{subsection-asymptotic-period} \

\ms Given $\G=(V,A)$, we will use the distance $d_{\max}(\x,\y):=\max_{i}|\x_i-\y_i|$ 
in both $[0,1]^V$ and $[0,1]^A$. As usual, 
$d_{\max}(\x, E):=\inf\{d_{\max}(\x,\y):\ \y\in E\}$ for $\x\in [0,1]^F$ and 
$E\subset [0,1]^F$. First notice that for $\G,\sigma$ and $a$ fixed, the set
\begin{equation}~\label{equation-open-gsa}
O_{\G,\sigma,a}:=
\left\{(T,\x)\in [0,1]^A\times [0,1]^V:\ 
\inf_{t\geq 0} d_{\max}\left(F^{t}_{\G,\sigma,T,a}(\x),\Delta_T\right) > 0 
\right\}
\end{equation}
is open. Indeed, if $(T,\x)\in O_{\G,\sigma,a}$, then 
$d_{\max}\left(F^{t}_{\G,\sigma,T,a}(\x),\Delta_T\right) > 3\epsilon$, for some 
$\epsilon > 0$ and all $t \geq 0$. The orbit does not change if we perturb $T$, 
i.~e., for each $T'\in B_{\epsilon}(T)$ we have $F^t_{\G,\sigma,T',a}(\x)
=F^t_{\G,\sigma,T,a}(\x)$ for all $t\in\nn$. Furthermore, since
$d_{\max}\left(F^t_{\G,\sigma,T',a}(\x),\Delta_{T'}\right) > 2\epsilon$, 
then $d_{\max}\left(F^{t}_{\G,\sigma,T',a}(\y),\Delta_{T'}\right) 
> \epsilon$ for all $\y\in B_{\epsilon}(\x)$, therefore 
$B_{\epsilon}(T)\times B_{\epsilon}(\x)\subset O_{\G,\sigma,a}$. 

\ms Since for each $\G,\sigma$ and $a$ fixed, the set $O_{\G,\sigma,a}$ is open 
in $[0,1]^A\times [0,1]^V$, then
\begin{equation}~\label{equation-open-av}
O_{a,V}:=\bigcup_{\G}\bigcup_{\sigma\in \{-1,1\}^A}\bigcup\left\{(F_{\G,\sigma,T,a},\x):\ 
(T,\x)\in O_{\G,\sigma,a}\right\}
\end{equation}
is measurable in $\F_{a,V}\times [0,1]^V$. 

\ms \subsubsection{The asymptotic period is measurable in 
$O_{a,V}$}~\label{subsubsection-is-measurable} \

\ms Fix $\G,\sigma$ and $a$, and take $(T,\x)\in O_{\G,\sigma,a}$ such that
$P(F_{\G,\sigma,T,a},\x)=\tau$. There exists $\hat\y\in \per_{\tau}(F_{\G,\sigma,T,a})$ and $N\in\nn$,
such that 
\[
d_{\max}\left(F^{N}_{\G,\sigma,T,a}(\x),\hat\y\right) < \epsilon := 
       \frac{\inf_{t\geq 0}d_{\max}\left(F^{t}_{\G,\sigma,T,a}(\x),\Delta_T\right)}{3}.
\]
Hence, for each $t\geq 0$ we have $d_{\max}\left(F^{N+t}_{\G,\sigma,T,a}(\x),
F^{t}_{\G,\sigma,T,a}(\hat\y)\right) <  a^t \epsilon$. By the same argument as in the
previous paragraph, 
$d_{\max}\left(F^{N+t}_{\G,\sigma,T',a}(\y), F^{t}_{\G,\sigma,T',a}(\hat\y)\right) 
< a^t\epsilon + a^{N+t}\epsilon < 2a^t\epsilon$ for all $t\geq 0$. 
Therefore $P(F_{\G,\sigma,T',a},\y) = \tau$
for all $(T',\y)\in B_{\epsilon}(T)\times B_{\epsilon}(\x)$. In this way we have
proved that for each $\G,\sigma$ and $a$ fixed, the set 
\[O^{(\tau)}_{\G,\sigma,a}:=\left\{(T,\x)\in [0,1]^A\times [0,1]^V:\ 
P(F_{\G,\sigma,T,a},\x) = \tau\right\}
\]
is an open set in $[0,1]^A\times [0,1]^V$, therefore
\[
P^{-1}\{\tau\}\cap O_{a,V}:=\bigcup_{\G}\bigcup_{\sigma\in \{-1,1\}^A}\bigcup
\left\{(F_{\G,\sigma,T,a},\x):\ (T,\x)\in O^{(\tau)}_{\G,\sigma,a}\right\}
\]
is measurable in $\F_{a,V}\times [0,1]^V$.

\ms \subsubsection{The asymptotic period is finite in 
$O_{a,V}$}~\label{subsubsection-is-finite} \

\ms For $(F_{\G,\sigma,T,a},\x)\in O_{a,V}$ let 
$\hat\x\in {\rm closure}\left(\{F^t_{\G,\sigma,T,a}(\x):\ t\in \nn\}_{t\in \nn}\right)$ 
be such that
\[
\inf_{t\geq 0}d_{\max}\left(F^t_{\G,\sigma,T,a}(\x),\Delta_T\right)=
d_{\max}\left(\hat\x,\Delta_T\right)=\epsilon >0.
\]
Notice that $(F_{\G,\sigma,T,a},\hat\x)\in O_{a,V}$ too, and 
$\inf_{t\geq 0}d_{\max}\left(F^t_{\G,\sigma,T,a}(\hat\x),\Delta_T\right)=
d_{\max}\left(\hat\x,\Delta_T\right)=\epsilon$.

\ms For $N\in\nn$ such that 
$d_{\max}\left(F^{N}_{\G,\sigma,T,a}(\x),\hat\x\right)<\epsilon$ we have
$d_{\max}\left(F^{N+t}_{\G,\sigma,T,a}(\x),F^{t}(\hat\x)\right)\leq a^t\epsilon$
for all $t\geq 0$. Now, for $\tau\in \nn$ such that $a^{\tau}<1/4$ and
$d_{\max}\left(F^{N+\tau}_{\G,\sigma,T,a}(\x),\hat\x\right) < \epsilon/4$, we have 
\begin{eqnarray*}
d_{\max}\left(F^{\tau}_{\G,\sigma,T,a}(\hat\x),\hat\x\right)
&\leq& d_{\max}\left(F^{N+\tau}_{\G,\sigma,T,a}(\x),\hat\x\right)+
     d_{\max}\left(F^{N+\tau}_{\G,\sigma,T,a}(\x),F^{\tau}(\hat\x)\right)\\
& < & \epsilon/4+ a^{\tau}\epsilon < \epsilon/2.
\end{eqnarray*}
Since $\inf_{t\geq 0}d_{\max}\left(F^t_{\G,\sigma,T,a}(\hat\x),\Delta_T\right)=\epsilon$, then 
\[
d_{\max}\left(F^{k\tau}_{\G,\sigma,T,a}(\hat\x),F^{(k-1)\tau}_{\G,\sigma,T,a}(\hat\x)\right)
< \frac{a^{(k-1)\tau}\epsilon}{2} < \frac{\epsilon}{2^{2k-1}},
\]
for all $k\in \nn$. Hence $\hat\y:=\lim_{k\to\infty}F^{k\tau}_{\G,\sigma,T,a}(\hat\x)$ 
exists and satisfies $F^{\tau}_{\G,\sigma,T,a}(\hat\y)=\hat\y$. Furthermore, since
$d_{\max}\left(\hat\x,\hat\y\right)< \epsilon\sum_{k=1}^{\infty}2^{1-2k}=2\epsilon/3$, then 
\[
d_{\max}\left(F^{N+k\tau}_{\G,\sigma,T,a}(\x),\hat\y\right) 
\leq 
d_{\max}\left(F^{N+k\tau}_{\G,\sigma,T,a}(\x),F^{k\tau}_{\G,\sigma,T,a}(\hat\x)\right)+
d_{\max}\left(F^{k\tau}_{\G,\sigma,T,a}(\hat\x),\hat\y\right)
< 
\frac{5a^k\epsilon}{3},
\]
for each $k\in\nn$, which implies 
$\lim_{t\to\infty}\left|F^{t}_{\G,\sigma,T,a}(\x)-F^{t}_{\G,\sigma,T,a}(\y)\right|=0$, 
where $\y=F^n_{\G,\sigma,T,a}(\hat\y)$ with $n\equiv -N\mod\tau$. It follows 
from here that $P(F_{\G,\sigma,a,T},\x)=\tau$. In addition, since
$\hat\y\in{\rm closure}\left(\{F^t_{\G,\sigma,T,a}(\x):\ t\in \nn\}_{t\in \nn}\right)$,
then 
$\min_{0\leq t < \tau}d_{\max}\left(F_{\G,\sigma,a,T}^t(\hat\y),\Delta_T\right)=\epsilon$.

\ms \subsubsection{The complement of $O_{a,V}$ has zero measure for
$a < \max_{u\in V}(\Id(u)+1)^{-1}$}~\label{subsubsection-density-argument}\

\ms Fix $\G,\sigma$, $a$, and for $(T,\x)\notin O_{\G,\sigma,a}$ let
$D^{(t)}:=D_{\G,\sigma,a}(F^{t}_{\G,\sigma,T,a}(\x))$ for each $t\in\nn$.  
Then, there exists a sequence $\{t_k\in \nn\}_{k\in \nn}$, and $u,v\in V$, such that 
\begin{eqnarray}~\label{equation-threshold-type}
T_{(u,v)} & = & \lim_{k\to\infty} \left(a^{t_k}\x_u+(1-a)
\sum_{t=1}^{t_k}a^{t-1}D_u^{(t_k-t)}\right) \nonumber\\ 
        & = & (1-a) \lim_{k\to\infty} \left(\sum_{t=1}^{t_k}a^{t-1}D_u^{(t_k-t)}\right). 
\end{eqnarray}

\ms For each $\x\in [0,1]^V$ fixed, define the fiber
\begin{equation}~\label{equation-the-fiber}
C_{\G,\sigma,a}(\x) :=  \left\{T\in [0,1]^A:\ (T,\x)\notin O_{\G,\sigma,a}\right\}.
\end{equation}
Then, according to Equation~\eqref{equation-threshold-type} we have
\[
C_{\G,\sigma,a}(\x)\subset \bigcup_{(u,v)\in A} 
\left\{T\in [0,1]^A:\ T_{(u,v)}\in \frac{1-a}{\Id(u)} 
\times {\rm closure} \left(\Omega_{\G,a,u}\right) \right\},  
\]
where
\[
\Omega_{\G,a,u}:=\bigcup_{N=1}^{\infty}
\left\{\sum_{s=0}^{N-1}a^{s}\k_{s}:\ \left(\k_s\right)_{s=0}^{N-1}\in 
\{0,1,\dots,\Id(u)\}^N\right\}.
\]

\ms If $a < (\Id(u)+1)^{-1}$ for some $u\in V$, then
${\rm closure}\left(\Omega_{\G,a,u}\right)$ is a Cantor set of dimension 
$\log(\Id(u)+1)/\log(1/a) < 1$. Hence, Fubini's Theorem implies,
for $a < \max_{v\in V}(\Id(u)+1)^{-1}$, that the fiber $C_{\G,\sigma,a}(\x)$ 
has zero $\#A$--dimensional volume for each $\x\in [0,1]^V$.

\ms \subsubsection{The complement of $O_{a,V}$ has zero measure for
$a \geq \max_{u\in V}(\Id(u)+1)^{-1}$}~\label{subsubsection-complexity-argument}\

\ms \ms The fact that $\vol \left( C_{\G,\sigma,a}(\x)\right) =0$ for 
each $\x\in [0,1]^A$ and arbitrary $a\in(0,1)$, 
derives from a result by Kruglikov and Rypdal. 
Theorem 2 in~\cite{Kruglikov&Rypdal} gives an upper bound for the topological entropy of the 
symbolic system associated to a non--degenerated piecewise affine map. The upper bound is
related to the exponential rate of angular expansion under the action of the map. 
For piecewise affine conformal maps, which is the case of discrete--time regulatory networks, this 
upper bound vanishes. Hence,  Kruglikov and Rypdal's theorem directly applies to our case, 
implying that
\begin{equation}~\label{equation-entropy}
          \limsup_{N\to\infty}\frac{\log \#\L_{(T,\x,u)}^{(N)}}{N}=0,
\end{equation}
where for each $\G,\sigma$, and $a$ fixed,
\[\L_{(T,\x,u)}^{(N)}:=
\left\{\left(\k_{t+s}\right)_{s=0}^{N-1}\in \{0,1,\ldots,\Id(u)\}^N:\
D_{\G,\sigma,T}(F^{t}_{\G,\sigma,T,a}(\x))_u=\frac{\k_{t}}{\Id(u)} \ \forall t\in\nn
\right\}. 
\]
Therefore, for each $T$, there exists $N=N(T,\x,u)\in \nn$ such that 
$\#\L_{(T,\x,u)}^{(N)}\leq ((1+a)/2a)^{N}$. Fix $u\in V$, and for each 
$N\in \nn$ and $\D\subset \{0,1,\ldots,\Id(u)\}^N$, define
\[
\T_{\D}:=\left\{x=\sum_{t=0}^{\infty}a^t\k_t:\ 
(\k_{(k+1)N-1},\ldots,\k_{kN})\in \D \ \forall k \geq 0\right\}.
\]
With this we have, for the fiber $C_{\G,\sigma,a}(\x)$ defined 
in~\eqref{equation-the-fiber}, the inclusion 
\begin{equation}~\label{equation-inclusion}
C_{\G,\sigma,a}(\x)\subset \bigcup_{(u,v)\in A}
                             \bigcup_{N\in\nn}
              \bigcup_{\D\subset \{0,1,\ldots,\Id(u)\}^N\atop \#\D\leq ((1+a)/2a)^{N} }
                   \left\{T\in [0,1]^A:\ T_{(u,v)}\in \frac{1-a}{\Id(u)}\times
                    {\rm closure}(\T_{\D})\right\}.
\end{equation}

\ms Now, for each $m\in \nn$, a prefix 
$\k\in \left(\L_{(T,\x,u)}^{(N)}\right)^m$ defines 
an interval
\[I_{\k}:=
\left\{x=\sum_{t=0}^\infty \omega_ta^t:\ 
(\omega_t)_{t=0}^{mN-1}=(\k_{mN-t})_{t=1}^{mN}\right\},
\]
and we have 
\[
{\rm closure}(\T_{\D})=\bigcap_{m=1}^{\infty} 
\left(\bigcup_{\k\in \D^m} {\rm closure}\left(I_{\k}\right)\right).
\]

\ms Since $\vol(I_{\k})=\vol({\rm closure}(I_{\k}))=\Id(u)a^{mN}/(1-a)$
for each $\k\in \left(\L_{(T,\x,u)}^{(N)}\right)^m$, then
\[
\vol({\rm closure}(\T_{\D})) \leq  (\#\D)^{m}a^{mN}\frac{ \Id(u)}{1-a},
\]
for each $m\in \nn$. Hence, $\vol({\rm closure}(\T_{\D}))=0$ for each
$\D\subset \{0,1,\ldots,\Id(u)\}^N$ such that $\#\D\leq ((1+a)/2a)^{N}$.
Taking into account~\eqref{equation-inclusion}, Fubini's Theorem implies, that the 
fiber $C_{\G,\sigma,a}(\x)$ has zero $\#A$--dimensional volume for each $\x\in [0,1]^V$.

\ms Using once again Fubini's, we finally conclude,
\[
\pp_{a,\eta}(\F_{a,V}\times [0,1]^V\setminus O_{a,V})\leq\sum_{\G\in\GG_V}
\sum_{\sigma\in\{-1,1\}^A}\pp_G(\G)\pp_{A,\eta}{\sigma}\int_{\x\in[0,1]^V} 
\vol\left(C_{\G,\sigma,a}(\x)\right) \ d\x=0.
\]

\endproof

\ms \begin{nota} In the statement of Theorem~\ref{theorem-asymptotic-period} we 
can replace the probability measure $\pp_{a,\eta}$ by any other Borel measure $\pp$
in $\F_{a,V}\times [0,1]^V$ such that, for each $\G\in\GG_V, \sigma\in \{-1,1\}^A$, 
and $\x \in [0,1]^V$ fixed, the $\#A$--dimensional projection 
\[
\pp_{T}(\J):=\pp\{(F_{\G,\sigma,T,a},\x)\in \F_{a,V}\times [0,1]^V:\ 
T\in \J\}
\]
is absolutely continuous with respect to the Lebesgue measure.
\end{nota} 

\ms \begin{nota} 
From the arguments developed in 
Paragraph~\ref{subsubsection-density-argument}, it follows that for 
$a<\min_{u\in V}(\Id(u)+1)^{-1}$, the set 
$O_{\G,\sigma,a}$ defined in~\eqref{equation-open-gsa} is open and dense. This 
ensures that orbits are generically uniformly far from the discontinuity set.
As a consequence, periodic attractors not intersecting the discontinuity set
are generic for $a<\min_{u\in V}(\Id(u)+1)^{-1}$. 
On the other hand, the argument in Paragraph~\ref{subsubsection-complexity-argument} 
ensures, as long as we consider distributions absolutely 
continuous with respect to Lebesgue for thresholds, that 
orbits are almost surely uniformly far from the discontinuity set.
Hence, in the general case, almost all initial conditions converge to
periodic attractors not intersecting the discontinuity set. 
\end{nota}

\bs \subsection{Proof of Theorem~\ref{theorem-modularity}}~\label{modularity}

\ms This theorem consist of two claims. The first one establishes the stability of 
the oscillatory subnetworks, while the second one gives the equivalence between the 
dynamics restricted to that subnetwork and the dynamics supported by the subnetwork 
considered as an isolated system. Both claims follow from direct computations.

\ms \subsubsection{Stability of the oscillatory 
subnetworks}~\label{subsubsection-stability}\

\ms Since $\pp_{a,\eta}(\G_\osc=\bG) >0$, then, there exists $\G\in\GG_V$,
$\sigma\in\{-1,1\}^A$, and $(\tT,\tx)\in O_{\G,\sigma,a}$, such that 
$\G_\osc(F_{\G,\sigma,\tT,a},\tx)=\bar\G$. 
As proved in Paragraph~\ref{subsubsection-is-finite}, there exists $0 <\epsilon < 1/2$, 
$\tau\in\nn$, and $\ty \in\per_\tau(F_{\G,\sigma,\tT,a})$, 
such that 
$\lim_{t\to\infty}\left|F^t_{\G,\sigma,\tT,a}(\tx)-
F_{\G,\bsigma,\tT,a}^t(\ty)\right|=0$ and
\[\inf_{t\geq 0}
d_{\max}\left(F_{\G,\sigma,\tT,a}^t(\tx),\Delta_{\tT}\right)=
\min_{0\leq t < \tau}d_{\max}
\left(F_{\G,\sigma,\tT,a}^t(\ty),\Delta_{\tT}\right)=2\epsilon.
\]
From this it readily follows that $\G_\osc(F_{\G,\sigma,T,a},\y)=\bar\G$
for each $(T,\y)\in B_\epsilon(\tT)\times B_\epsilon(\ty)$.

\ms Let us now extend $\G$, $\sigma$, and $\tT$, and redefine $\ty$ at each vertex 
$v\in \bV$ such that $\Idosc(v)=\Id(v)$ (remind that $\Idosc(v):=\#\{u\in V:\ (u,v)\in \bA\equiv A_\osc\}$). 
Choose $u\notin \bV$, and include the arrow $(u,v)$ in $A$. Then, if 
$\ty_u < 1-2\epsilon$, define $\sigma_{(u,v)}=1$ and $\tT_{(u,v)}=1-\epsilon$.
Otherwise, if $\ty_u > 2\epsilon$, make $\sigma_{(u,v)}=-1$ and $\tT_{(u,v)}=\epsilon$.
Finally, redefine $\ty_v=\Id(v)\ \by_v/(\Id(v)+1)+1/(\Id(v)+1)$ at that vertex.
A simple computation shows that, in the redefined system, $\G_\osc(F_{\G,\sigma,T,a},\y)=\bG$, 
for each $(T,\y)\in B_\epsilon(\tT)\times B_\epsilon(\ty)$. Therefore, taking into account
that by hypothesis $\pp_G(\G)>0$ for all $\G\in\GG_V$, we can assume without lost
of generality that $\Idosc(v) < \Id(v)$ for all $v\in \bV$.

\ms By definition, if $(u,v)\in A\setminus \bA$ then 
$\theta_{(u,v)}:=H(\bsigma_{(u,v)}(\ty^t_u-\tT_{(u,v)}))$ remains constant in time, as well as
$D(v):=\sum_{(v,u)\in A\setminus \bA} \theta_{(v,u)}$ for each $v\in V$. 
Let us now modify $\sigma$ and $\bT$ for the arrows leaving $\bG$. 
For $u\in \bV$, and each
$v\notin \bV$, if $D(u)>0$, make $\tT_{(u,v)}\mapsto \tT_{(u,v)}\times D(u)/\Id(u)$, and
\[
\sigma_{(u,v)}\mapsto \left\{\begin{array}{ll} \sigma_{(u,v)} 
                             & \text{ if } \theta_{(u,v)}+\sigma_{(u,v)}\in \{-1,2\} \\
                                               -\sigma_{(u,v)}
                              & \text{ if } \theta_{(u,v)}+\sigma_{(u,v)}\in\{0,1\}.\end{array}
                         \right.
\]
Otherwise, if $D(u)=0$, then 
$\tT_{(u,v)}\mapsto \tT_{(u,v)}\times (\Id(u)-\Idosc(u))/\Id(u)+\Idosc(u)/\Id(u)$, and
\[
\sigma_{(u,v)}\mapsto \left\{\begin{array}{ll} -\sigma_{(u,v)} 
                             & \text{ if } \theta_{(u,v)}+\sigma_{(u,v)}\in \{-1,2\} \\
                                               \sigma_{(u,v)}
                              & \text{ if } \theta_{(u,v)}+\sigma_{(u,v)}\in\{0,1\}.\end{array}
                         \right.
\]
This modification in $\sigma$ and $\tT$ uncouples the dynamics on $\bV$ from that on 
$V\setminus \bV$, so that now we can set $\sigma|_{\bA}\equiv\bsigma$ for arbitrary 
$\bsigma\in\{-1,1\}^{\bA}$. In this way we obtain a digraph extension $\bG_\ext:=\G\supset\bG$, and
for each $\bsigma\in\{-1,1\}^{\bA}$ an extension $\bsigma_\ext:=\sigma\in \{-1,1\}^A$ with
$\bsigma_\ext|_{\bA}=\bsigma$. For any of this extensions, a direct computation shows that
$\lim_{t\to\infty}F_{\bG_\ext,\bsigma_\ext,T,a}^t(\x)_v=D(v)/\Id(v)=\by_v$ for 
$v\notin \bV$, and 
$F^t_{\bG_\ext,\bsigma_\ext,T,a}(\x)_v\in\left[D(v)/\Id(v),(D(v)+\Idosc(v))/\Id(v)\right]$, 
for all $t\geq 0$ and $v\in \bV$, whenever
\[\x\in \bJ\times \J':=
\left(\prod_{v\notin \bV} [\ty_v-\epsilon,\ty+\epsilon]\cap [0,1]\right)
\times
\left(\prod_{v\in \bV}\left[\frac{D(v)}{\Id(v)},\frac{D(v)+\Idosc(v)}{\Id(v)}\right]\right),
\] and
\[T\in \bI\times\I':=\left(\prod_{(u,v)\in A\setminus \bA}[\tT_{(u,v)}-\epsilon,\tT_{(u,v)}+\epsilon]\cap [0,1]\right)
\times \left(\prod_{(u,v)\in \bA} \left[\frac{D(u)}{\Id(u)},\frac{D(u)+\Idosc(u)}{\Id(u)}\right]\right).
\]
Therefore, $\G_\osc(F_{\bG_\ext,\bsigma_\ext,T,a},\x)\subseteq \bG$ for each 
$(T,\x)\in \I\times\J:=\left( \bI\times\I' \right)\times \left( \bJ\times \J' \right)$,
and $\bsigma\in\{-1,1\}^{\bA}$.

\ms \subsubsection{Equivalence to the subnetwork considered as an isolated system}\
 
\ms We will use the same notation as in Subsubsection~\ref{subsubsection-stability}. 
Let $\sigma:=\bsigma_\ext$, and each $(T,\x)\in \left(\bI\times\I'\right)\times 
\left(\bJ\times\J'\right)$, and $t\in \nn$, let 
$\x^t:=F_{\bG_\ext,\bsigma_\ext,T,a}^t(\x)$. For each $v\in V$, the term
\[
D(v):=\sum_{(u,v)\in A\setminus \bA} 
H\left(\sigma_{(u,v)}\left(\x^t_u-T_{(u,v)}\right)\right)=
\sum_{(u,v)\in A\setminus \bA} 
H\left(\sigma_{(u,v)}\left(F_{\G,\sigma,\tT,a}^t(\ty)_u-\tT_{(u,v)}\right)\right),
\]
remains constant in time. Since $\sigma|_{\bA}\equiv \bsigma$, 
\[
\x^{t+1}_v=a \x^t_v+
\frac{\Idosc(v)}{\Id(v)}\left(\frac{(1-a)}{\Idosc(v)} 
\sum_{(u,v)\in \bA}H(\bsigma_{(u,v)}(\x^t_u-T_{(u,v)}))\right)+
(1-a)\frac{D(v)}{\Id(v)},
\]
which can be rewritten as 
\begin{equation}~\label{equation-equivalence}
\Phi_V(\x^{t+1})_v=a \Phi_V(\x^t)_v+\frac{1-a}{\Idosc(v)}\sum_{(u,v)\in \bA}
H\left(\bsigma_{(u,v)}\left(\Phi_V(\x^t)_u-\Phi_A(T)_{(u,v)}\right)\right),
\end{equation}
where $\Phi_A:\bI\times\I'\to [0,1]^{\bA}$ and  $\Phi_V:\bJ\times\J'\to [0,1]^{\bV}$  
are affine transformations defined as follows. We have
$\Phi_A(\bT\times T')=\D_{\bA}\bT+\C_{\bA}$, with
\[
(\D_{\bA}\bT)_{(u,v)}=\frac{\Id(u)}{\Idosc(u)}\times \bT_{(u,v)}, 
\ \text{ and }\ (\C_{\bA})_{(u,v)}=-\frac{D(u)}{\Idosc(u)}, 
\] for each $(u,v)\in \bA$, and  
$\Phi_V(\bx\times \x')=\D_{\bV}\bx+\C_{\bV}$, with
\[
(\D_{\bA}\bx)_u=\frac{\Id(v)}{\Idosc(v)}\times \bx_v,
\ \text{ and }\ (\C_{\bV})_v=-\frac{D(v)}{\Idosc(v)}, 
\] for each $v\in\bV$. The result follows from Equation~\eqref{equation-equivalence}, 
which establishes 
\[
F_{\bG,\bsigma,\Phi_A(T),a}(\Phi_V(\x))=\Phi_V(F_{\bG_\ext,\bsigma_\ext,T,a}(\x))
\] 
for all $(T,\x)\in \I\times \J$, and from the fact that $\Phi_A$ and $\Phi_V$ are affine 
and surjective.

\endproof

\ms \begin{nota} In the statement of Theorem~\ref{theorem-modularity}, instead of
the hypothesis $\pp_G(\G) >0$ for all $\G\in\GG_V$, we could have directly 
supposed that the oscillatory subnetwork $\bG:=(\bV,\bA)$ is such that $\Idosc(v) < \Id(v)$ 
for each $v\in \bV$. By assuming this, we avoid the condition $\bV\subsetneq V$ as well.
This alternative formulation allows us to consider the emergence of modularity in 
Barab\'asi--Albert random networks, for which $\pp_{\rm BA}(\G)=0$ for some
digraphs $\G\in\GG_V$.
\end{nota}

\bs \section{Final Comments and Conclusions} \

\ms From our point of view, one of the most important theoretical problem 
in regulatory dynamics concerns the relations between the structure of the underlying 
network and the possible dynamical behaviors the system generates. We have 
already established bounds relating the structure of the underlying network, 
and the growth of distinguishable orbits~\cite{Lima&Ugalde06}. Here, 
following~\cite{Volchenkov&Lima05}, we have considered regulatory networks 
with interactions and initial conditions chosen at random at the beginning 
of the evolution. Within this approach, individual orbits become 
elements in a sample space, and the characteristics of the asymptotic
behavior can be considered as orbit dependent random variables. The structure of
the underlying network, random in this approach, is encoded in a 
probability distribution over a set of directed graphs. In two interesting 
cases we have explored, the asymptotic oscillations concentrates on a relatively small 
subnetwork, whose structure depends on both the proportion of inhibitory interactions, 
and the contraction rate. We have proved that 
the dynamics observed on this subnetwork is equivalent to the dynamics supported 
by the subnetwork considered as an isolated system. 
We interpret this as the emergence of modularity. 
Modularity allows us to predict asymptotic periods in regulatory networks 
admitting disconnected oscillatory subnetworks. 

\ms We want to remark Theorem~\ref{theorem-asymptotic-period}, which has important 
technical implications. It ensures, as long as we consider distributions absolutely 
continuous with respect to Lebesgue for thresholds and initial conditions, that a 
random orbit converges to a periodic attractor. Furthermore, this periodic attractor 
does not intersect the discontinuity set. We have also proved, in the case of small 
contraction rates, that the periodic attractors are generic. This can be deduced 
from the argument in Paragraph~\ref{subsubsection-density-argument}.

\ms As mentioned above, in order to determine the distribution of the asymptotic period, 
and the characteristics of the oscillatory subnetwork, a detailed numerical study
of the examples presented in Section~\ref{section-examples} is required. Heuristic 
computations could guide those numerical studies, and it is our intention to proceed 
in this direction in a subsequent work.

\bs 

\end{document}